\documentclass[11pt,twoside]{article} 
\usepackage{amsmath}
\usepackage{latexsym} 
\usepackage{tikz} 
\input{amssym.def} 
\input{amssym} 
\usepackage{color}
\usepackage{multirow}
\newcommand\blue{\color{blue}}
\newcommand\black{\color{black}}

\newfont{\fra}{eufm10 scaled 1095} 
\newfont{\Bb}{msbm10 scaled 1095} 
\newfont{\Bbg}{msbm10 scaled 1680} 
\newcommand\CC{{\mbox{\Bb C}}} 
\newcommand\RR{{\mbox{\Bb R}}} 
\newcommand\NN{{\mbox{\Bb N}}} 
\newcommand\ZZ{{\mbox{\Bb Z}}} 
\newcommand\SSS{{\mbox{\Bb S}}} 
\newcommand\Z{{\Bbb Z}}
\newcommand\QQ{{\mbox{\Bb Q}}}

\newcommand\fg{{\frak{g}}} 
 
\newcommand\fri{{\frak i}}

\newcommand\fn{{\frak n}}

\newcommand\fq{{\frak q}}

\newcommand\fz{{\frak z}} 

\newcommand\cO{{{\cal O}}} 
\newcommand\cS{{\cal S}}
\newcommand\cL{{\cal L}} 
 
\newcommand\cB{{\cal B}} 
\newcommand\cP{{\cal P}}
\newcommand\cR{{\cal R}}
\newcommand\cG{{\cal G}}
\newcommand\cK{{\cal K}}
\newcommand\bB{{\bf B}}
\newcommand\bH{{\bf H}}
\newcommand\bz{{\bf z}}
\newcommand\ph{\varphi} 
 
\newcommand\eps{\varepsilon} 

\newcommand\fosc{\frak o \frak s \frak c_{1,0}}

\newcommand{\GL}{\mathop{{\rm GL}}} 
\newcommand{\SL}{\mathop{{\rm SL}}} 
\newcommand{\SO}{\mathop{{\rm SO}}}

\newcommand{\Aut}{\mathop{{\rm Aut}}} 
\newcommand{\Osc}{{{\rm Osc_{1,0}}}} 
\newcommand{\dOsc}{{\rm Osc}^r}

\newcommand{\Id}{{{\rm id}}} 
\newcommand{\ad}{{{\rm ad}}} 
 
\newcommand{\tr}{\mathop{{\rm tr}}} 
\newcommand{\Ad}{{{\rm Ad}}}

\newcommand{\diag}{\mathop{{\rm diag}}}

\newcommand{\Span}{{{\rm span}}} 
\newcommand\ip{{\langle\cdot \,,\cdot \rangle}}

\newcommand\la{{\langle}} 
\newcommand\ra{{\rangle}} 
\newcommand\proof{{\sl Proof. }} 
\newcommand{\qed}{\hspace*{\fill}\hbox{$\Box$}\vspace{2ex}} 
\newcommand{\qedohne}{\hspace*{\fill}\hbox{$\Box$}} 
\newtheorem{theo}{Theorem}[section] 
\newtheorem{pr}[theo]{Proposition}
\newtheorem{de}[theo]{Definition}
\newtheorem{ex}[theo]{Example}
\newtheorem{re}[theo]{Remark}
\newtheorem{co}[theo]{Corollary}
\newtheorem{lm}[theo]{Lemma}
\newcommand{\benur}{\begin{enumerate}[label=(\roman*)]}
\begin{document} 
\title{Lattices in the four-dimensional split oscillator group} 
\author{Blandine Galiay and Ines Kath}
\maketitle 
\begin{abstract}
Besides the oscillator group, there is another four-dimensional non-abelian solvable Lie group that admits a bi-invariant pseudo-Riemannian metric. It is called split oscillator group (sometimes also hyperbolic oscillator group or Boidol's group). 
We parametrise the set of lattices in this group and develop a method to classify these lattices up automorphisms of the ambient group. We show that their commensurability classes are in bijection with the set of real quadratic fields.
\end{abstract}
\section{Introduction}
In this paper, we study lattices in a certain solvable four-dimensional Lie group. A lattice in a Lie group $G$ is a discrete subgroup $L$ with the property that the quotient space $L\backslash G$ has finite invariant measure. Lattices in a solvable Lie group are uniform, that is, the quotient space is compact. We are particularly interested in the case where $G$ is solvable and admits a bi-invariant metric. This is motivated by the following
result by Baues and Globke \cite{BG}: Let $M$ be a compact pseudo-Riemannian manifold, and let $G$ be a connected solvable Lie group of isometries acting transitively on $M$. Then $M$ is diffeomorphic to a quotient of $G$ by a lattice and the pseudo-Riemannian metric on $M$ pulls back to a bi-invariant metric on $G$. 

The smallest non-abelian solvable Lie groups that admit a bi-invariant metric are four-dimensional. One of them is the  (ordinary) oscillator group ${\rm Osc_{0,1}}$, which admits a Lorentzian bi-invariant metric. Recall that the oscillator group is a semi-direct product of the Heisenberg group $H$ and the real line, where $\RR$ acts trivially on the centre $Z(H)$ of $H$ and by rotations on $H/Z(H)$. The other group has a similar structure. It is also a semi-direct product of the Heisenberg group and the real line, but now $\RR$ acts by hyperbolic rotations. It is called {\it split oscillator group}. Other names in the literature are {\it hyperbolic oscillator group} or {\it Boidol's group}. We will denote it by $\Osc$. The bi-invariant metric on this group has signature $(2,2)$.

Lattices in the oscillator group ${\rm Osc}_{0,1}$ were classified up to automorphisms of ${\rm Osc}_{0,1}$ in~\cite{F} and up to inner automorphisms of ${\rm Osc}_{0,1}$ in~\cite{FK}. Lattices in ${\rm Osc}_{0,1}$ are not rigid, that is, there occur continuous families of lattices that are isomorphic as discrete groups but cannot be mapped into each other by an automorphism of ${\rm Osc}_{0,1}$. The classification is given in form of a list. For each isomorphism class of the possible discrete groups, a bijective parametrisation of the set of lattices (modulo automorphisms of ${\rm Osc}_{0,1}$) is given whose underlying discrete group belongs to this class. 

The situation changes if we consider the split oscillator group. On the one hand it becomes easier since, unlike the ordinary oscillator group, the split one is completely solvable. Hence we can apply a general  rigidity result by Sa\^ito~\cite{S}, which yields that isomorphism classes of lattices modulo $\Aut(\Osc)$ are the same as the isomorphism classes of the underlying discrete subgroups. On the other hand, it is more complicated. While the classification of lattices in the ordinary oscillator group is linked to conjugacy classes of matrices of finite order in $\GL(2,\ZZ)$, now a description of conjugacy classes of matrices of infinite order is needed and such a description cannot be given in form of an explicit list. Thus we cannot expect to obtain a list of all lattices here. However, there are several classical methods to describe conjugacy classes in $\GL(2,\ZZ)$. Here we want to apply the one developed by Aicardi~\cite{A}.

A discrete group that can be embedded into $\Osc$ as a lattice will be called a discrete split oscillator group. Such a group is a semi-direct product of a discrete Heisenberg group $H_{\rm d}$ with an infinite cyclic group generated by an element $\delta$ with respect to a homomorphism $\phi:\ZZ\cdot\delta\rightarrow \Aut(H_{\rm d})$ for which the map $B:=\overline{\phi(\delta)}$ induced by $\phi(\delta)$ on $H_{\rm d} / Z(H_{\rm d} )\cong \ZZ^2$ is conjugate in $\GL(2, \RR)$ to a hyperbolic rotation.  

In this paper, 
\begin{itemize}
\item[-] we give a description of the set of isomorphism classes of discrete split oscillator groups, see Proposition~\ref{pre};
\item[-] we parametrise the set of all lattices in the split oscillator group, see Proposition~\ref{Ppar};
\item[-] for each lattice, we determine the isomorphism class of the underlying discrete split oscillator group, see Corollary~\ref{Bl};
\item[-]  using Proposition~\ref{pre} and the description of conjugacy classes in $\GL(2,\ZZ)$ developed in~\cite{A}, we derive a classification method for isomorphism classes of discrete split oscillator groups, see Theorem~\ref{theo};
\item[-] we apply this general method to obtain an explicit classification of isomorphisms classes of discrete split oscillator groups with  $T:=\tr B\le 7$, see the table in Section~\ref{S83};
\item[-] we give a bijection from the set of commensurability classes of discrete split
oscillator groups to the set of real quadratic fields, see Theorem~\ref{Tcomm}.
\end{itemize}

In the last section we relate our results to those obtained by Maeta  \cite{Ma}, who considers compact Clifford-Klein forms of solvable pseudo-Riemannian symmetric spaces. More exactly, Maeta considers compact quotients of solvable symmetric spaces by discrete subgroups of their transvection group. He proved that for signature $(2,2)$ such quotients exist only for two solvable symmetric spaces. One of these spaces is the split oscillator group endowed with its bi-invariant metric. On the other hand, each lattice of $\Osc$ gives rise to a compact quotient. However, we will see that these quotients cannot be obtained as a quotient of the symmetric space $X=\Osc$ by a discrete subgroup of its transvection group. Indeed, the transvection group $\hat G$ of $X$ is properly contained in the identity component $G$ of the isometry group. In particular, we will see that the subgroup of $G$ that consists of left-translations by group elements of $\Osc$ is not contained in $\hat G$. Consequently, we have to consider $\Osc$ as a homogeneous space of $G$ in order to understand $L\backslash \Osc$ as a Clifford-Klein form. In particular, the compact quotients of the symmetric space $\Osc$ considered in~\cite{Ma} are not of the form $L\backslash \Osc$. We will see that the quotients of $X$ constructed in~\cite{Ma} are straight in the sense of \cite{KO} and we will prove that there exist also non-straight  quotients. 
\section{The split oscillator group}  
\subsection{Definition and basic properties} \label{S21}
The split oscillator group is a semi-direct product of the Heisenberg group $H$ and the real line, where $\RR$ acts trivially on the centre $Z(H)$ of $H$ and by hyperbolic rotations on $H/Z(H)$. To make this more precise, we start by recalling the definition of the three-dimensional Heisenberg group.
\begin{de}
The Heisenberg group $H$ is defined as the group $\RR \times \RR^2$ with multiplication given by 
$$(z, \xi )\cdot(z', \xi ' ) = (z+z'+ \textstyle \frac12\omega(\xi, \xi '), \xi +\xi '),$$
where $\omega$ is the symplectic form defined by  $\omega (\xi , \xi ') = xy'-x'y$ for $\xi=(x,y)$, $\xi'=(x',y')$. 
\end{de}

The automorphism group of the Heisenberg group is isomorphic to $\GL(2,\RR)\rtimes \RR^2$. The elements $P\in\GL(2,\RR)$ and $\eta\in \RR^2$ act by the automorphisms  $\bar F_P$ and $\bar F_\eta$ given by 
\begin{eqnarray*}
\bar F_P (z,\xi) &=& (\det(P)\cdot z, P\xi),\\
\bar F_\eta (z,\xi) &=& \Ad_{(0,\eta)}(z,\xi)\ =\ (z+\omega(\eta,\xi),\xi)\,.
\end{eqnarray*}

In this paper, let $A$ be the matrix
$$A = \begin{pmatrix} 0 & 1 \\ 1 & 0 \end{pmatrix}\in{\SO}^+(1,1).$$ 
\begin{de}
The split oscillator group $\Osc$ is defined as the semi-direct product $H\rtimes \RR$, where $t\in \RR$ acts on $H$ by 
\begin{equation}\label{action}
 t.(z, \xi) = \bar F_{e^{tA}}(z,\xi)=(z, e^{tA} \xi ).  
\end{equation}
\end{de}
Hence the multiplication in $\Osc$ is given by
$$(z, \xi , t )\cdot(z', \xi ', t' ) = \left( z+z'+ \textstyle\frac12\omega(\xi, e^{tA} \xi '), \xi + e^{tA} \xi ', t+t' \right).$$
The group $\Osc$ admits a bi-invariant metric of signature (2,2). The induced $\ad$-invariant inner product $\ip$ on the Lie algebra $\fosc$ of $\Osc$ is given as follows. Since $\Osc=\RR\times\RR^2\times\RR$ as a set, we can identify $\fosc$ with $\RR\oplus\RR^2\oplus\RR$. Then
$$\big\langle (Z,X,T), (Z',X',T')\big\rangle =\la X,X'\ra_0 +ZT'+Z'T,$$
where $\ip_0= \diag(-1,1)$ is the standard pseudo-Euclidean scalar product of signature $(1,1)$ on $\RR^2$.

\subsection{The automorphism group of $\Osc$} \label{S22}
In this subsection we determine the automorphism group of the split oscillator group $\Osc$. The automorphism group of its Lie algebra was already described in \cite{Ti}. Of course, one could use this result to compute the automorphism group of $\Osc$. However, we prefer a more concise presentation similar to that for the ordinary oscillator group of arbitrary dimension given by Fi\-scher~\cite{F}. 

We define the group
$$\cS:=\{ S\in\GL(2,\RR)\mid SA=\mu AS,\ \mu=\pm1\}.$$
Then $\cS=\diag(1,-1)\cdot\cS^+=\cS^+\cdot\diag(1,-1)$ for $\cS^+:=\{ S\in\GL(2,\RR)\mid SA=AS\}$. 

Analogously to the case of the (ordinary) oscillator group one checks that the following maps are automorphisms of $\Osc$\,:
\begin{itemize}
\item[(i)] $F_u:\Osc\to\Osc$, $(z,\xi, t)\mapsto (z+ut,\xi,t)$ for $u\in\RR$, 
\item[(ii)] $F_\eta:\Osc\to\Osc$ for $\eta\in\RR^2$, where $F_\eta$ restricted to $H$ equals $\bar F_\eta$ and $$F_\eta(0,0,t)=\big(-\textstyle\frac12\omega(\eta,e^{tA}\eta),\, \eta-e^{tA}\eta, \,t\big),$$
\item[(iii)] $F_S:\Osc\to\Osc$, $(z,\xi, t)\mapsto (\det(S)z,S\xi,\mu t)$ for $S\in\cS$ such that $SA=\mu AS$.
\end{itemize}
\begin{pr} The automorphism group of $\Osc$ is isomorphic to $(\RR\times\RR^2)\rtimes \cS$, where $S\in \cS$ acts on $\RR\times\RR^2$ by $S\cdot (u,\eta)=(\mu \det (S) u, S\eta)$. The isomorphism is given by
\begin{eqnarray*}
(\RR\times\RR^2)\rtimes \cS &\longrightarrow & \Aut(\Osc) \\
(u,\eta,S)&\longmapsto & F_u\circ F_\eta \circ F_S. 
\end{eqnarray*}
\end{pr}
\proof Let $F$ be an automorphism of $\Osc$. The Heisenberg group $H$ is the commutator subgroup of $\Osc$. Hence $F$ leaves invariant $H$ and induces an automorphism of $\Osc/H$. In particular, $F(0,0,t)=(z(t),\xi(t),\mu\cdot t)$ for some $\mu \in\RR^*$. The condition $F(0,0,s)\cdot F(0,0,t)=F(0,0,s+t)$ yields
$$\xi(s+t)=\xi(s) + e^{\mu s A}\xi(t).$$
Differentiation with respect to $s$ gives 
\[\xi'(t)=\xi'(0)+\mu A(\xi(t))\]
in $s=0$. We define $\eta$ by $-\xi'(0)=\mu A\eta$ and obtain
\[\xi(t)=\eta-e^{\mu tA}\eta\]
because of $\xi(0)=0$. Now we consider $F_1:=F_{-\eta}\circ F$. We have
$$F_1(0,0,t)=F_{-\eta}(F(0,0,t))=F_{-\eta}(z(t),\eta-e^{\mu tA}\eta,\mu t)=(z_1(t),0,\mu t)$$
for some $z_1=z_1(t)$. Since $F_1$ is an automorphism, $z_1$ is linear, i.e., $z_1(t)=ct$ for some $c\in\RR$. Thus $F_2:=F_{-u}\circ F_1$ for $u:=c/\mu$ maps $(0,0,t)$ to $(0,0,\mu t)$. Moreover, $\bar F_2:=F_2|_H$ is an automorphism of the Heisenberg group. This implies $\bar F_2=\bar F_{\eta'}\circ \bar F_S$ for some $S\in\GL(2,\RR)$ and $\eta'\in\RR^2$. Since $F_2$ is an automorphism of $\Osc$, Eq.~(\ref{action}) implies $\bar F_2 \circ\bar F_{e^{tA}}=\bar F_{e^{\mu tA}}\circ \bar F_2$, hence
\[ \bar F_{\eta'}\circ \bar F_S\circ\bar F_{e^{t A}} = \bar F_{e^{\mu t A}}\circ \bar F_{\eta'}\circ \bar F_S= \bar F_{e^{\mu t A}\eta'}\circ \bar F_{e^{\mu t A}}\circ\bar F_S,\]
which gives $\eta'=e^{\mu t A}\eta'$ and $S e^{tA}=e^{\mu t A} S$ for all $t\in\RR$. Differentiating with respect to $t$, we get $\eta'=0$ and $SA=\mu AS$. In particular, $\mu=\pm1$. We obtain $$F_2(z,\xi,t)=(\bar F_2(z,\xi),\mu t)=(\bar F_S(z,\xi),\mu t)=F_S(z,\xi,t).$$
Finally, $F=F_{\eta}\circ F_1=F_{\eta}\circ F_{u}\circ F_2=F_{\eta}\circ F_{u}\circ F_S.$ \qedohne
\section{Basic notions and first examples of lattices}
\subsection{Lattices in the Heisenberg group}
Lattices in the $(2n+1)$-dimensional Heisenberg group were classified by Tolimieri~\cite{T}. However, note that, in \cite{T}, it is assumed that a lattice $\Gamma$ satisfies the additional technical condition $[\Gamma,\Gamma]=\Gamma\cap Z$, where $Z$ denotes the centre of the Heisenberg group. 

Let $e_1,e_2$ denote the standard basis of $\RR^2$. For each $r\in\NN_{>0}$, the elements $(1/r, 0)$, $(0,e_1)$ and $(0,e_2)$ generate a lattice $\Gamma_r$ in $H$. It is called standard lattice. By definition, we have
\begin{equation}\label{Gammar}
\Gamma_r=\left\{ (z,\xi)\in H \ \left| \  \xi=(\xi_1,\xi_2)\in\ZZ^2,\ z\in \textstyle \frac12\xi_1\xi_2+\frac1r \ZZ\right.\right\}.
\end{equation}
If $r$ is even, then this simplifies to $\Gamma_r=\left(\frac1r \ZZ\right)\times\ZZ^2$ (as sets).

An arbitrary lattice of $H$ equals $\bar F_\eta\bar F_P (\Gamma_r)$ for some $r\in\NN_{>0}$, $P\in\GL(2,\RR)$ and $\eta\in\RR^2$.

\begin{de}\label{df}
For each $r \in \NN_{>0}$, we define $\mathcal{B}_r$ to be the subgroup of $\RR^2\rtimes \GL(2,\RR)\cong\Aut(H)$ consisting of all elements $(\xi,B)$ that stabilise the lattice $\Gamma_r$, i.e., for which 
\begin{equation}\label{Br}
\bar F_\xi\bar F_B(\Gamma_r) =\Gamma_r
\end{equation}
holds.
\end{de}
\begin{lm}\label{LBr} Let $(\xi,B)$ be in $\RR^2\rtimes \GL(2,\RR)$, where $B=$\begin{small} $ \begin{pmatrix}a&b\\c&d \end{pmatrix}$\end{small}. Then $(\xi,B)$ is in $\cB_r$ if and only if $B$ is in $\GL(2,\ZZ)$ and the following condition for $\xi$ holds.
If $r$ is even, then $\xi \in \frac1r \ZZ \times \frac1r \ZZ $ and, if
$r$ is odd, then
\begin{equation}\label{EBr}
\xi\in \left\{   \begin{array}{ll} \frac1r \ZZ \times \frac1r \ZZ, & \mbox{if $ab$ and $cd$ are even}, \\[1ex]
\frac1r \ZZ \times (\frac{1}{2}+ \frac1r \ZZ),  & \mbox{if $ab$ is even and $cd$ is odd},\\[1ex]
(\frac{1}{2}+ \frac1r \ZZ) \times \frac1r \ZZ, & \mbox{if $ab$ is odd and $cd$ is even}.
\end{array}\right.
\end{equation}
\end{lm}
\proof We have 
\begin{eqnarray*}
\bar F_\xi\bar F_B\big((0,e_1)\big)&=&\big(\xi_1c-\xi_2a,(a,c)^\top\big) ,\\
\bar F_\xi\bar F_B\big((0,e_2)\big)&=&\big(\xi_1d-\xi_2b,(b,d)^\top\big).
\end{eqnarray*}
Suppose  that $\bar F_\xi\bar F_B(\Gamma_r) =\Gamma_r$. Then these equations imply that $a,b,c$ and $d$ are integers and we obtain
\begin{eqnarray*}
\bar F_\xi\bar F_B\big((0,e_1)\big)&=&(\xi_1c-\xi_2a-\textstyle \frac12 ac,0)(0,e_1)^a(0,e_2)^c,\\
 \bar F_\xi\bar F_B\big((0,e_2)\big)&=& (\xi_1d-\xi_2b-\textstyle \frac12 bd,0)(0,e_1)^b(0,e_2)^d.
 \end{eqnarray*}
In particular, (\ref{Gammar}) gives
\begin{equation}\label{odd}
\xi_1c-\xi_2a-\textstyle\frac12ac\in\frac1r\ZZ,\quad 
\xi_1d-\xi_2b-\textstyle\frac12bd\in\frac1r\ZZ.
\end{equation}
Conversely, if $B$ is an integer matrix and if (\ref{odd}) holds, then $\bar F_\xi\bar F_B(\Gamma_r) =\Gamma_r$. 
If $r$ is even, Eq.\,(\ref{odd}) is equivalent to $\xi \in \frac1r \ZZ \times \frac1r \ZZ$. If $r$ is odd, then (\ref{odd}) is equivalent to (\ref{EBr}).
\qedohne

\subsection{Examples of lattices in $\Osc$}
Let $B\in\SL(2,\ZZ)$ be an integer matrix with trace $T>2$. Then $B$ is conjugate to $e^{sA}$ for a unique $s\in\RR_{>0}$. Choose $P\in\GL(2,\RR)$ such that $B=P^{-1}e^{sA}P$ and let $\xi\in\RR^2$ be such that $(\xi,B)$ is in $\cB_r$, i.e., such that $\xi$ satisfies (\ref{EBr}). 

\begin{ex}\label{LxiP}
The subgroup $L_{r,s}(P,\xi)$ of $\Osc$ that is generated by the lattice $\bar F_{P}(\Gamma_r)$ of $H\subset\Osc$ and the element $\delta:=(0,P\xi,s)\in \Osc$ is a lattice in $\Osc$. 
\end{ex}
Indeed, let $F$ denote the conjugation by $\delta$. Then $F(H)=H$ and $F|_H=\bar F_{P\xi}\bar F_{e^{sA}}$\,. Hence 
$$
F\bar F_{P}(\Gamma_r)=\bar F_{P\xi}\bar F_{e^{sA}}\bar F_{P}(\Gamma_r)= \bar F_{P\xi}\bar F_{P}\bar F_{B}(\Gamma_r) = \bar F_{P} \bar F_\xi\bar F_{B}(\Gamma_r)= \bar F_{P}(\Gamma_r)$$
by Lemma~\ref{LBr}, which proves the claim.

In the next subsection we will define isomorphy of lattices. Later on, we will see that each lattice in $\Osc$ is isomorphic to one of the lattices $L_{r,s}(P,\xi)$.

\section{Isomorphy of lattices and rigidity}\label{S4}

Two lattices $L,L'$ in a group $G$ are called isomorphic if there exists an automorphism $F$ of $G$ such that $F(L)=L'$. If two lattices are isomorphic, then their underlying discrete groups are isomorphic. The converse is not true in general. In particular, it is not true for the (ordinary) oscillator group (\cite{F}, see also \cite{FK}). However, the split oscillator group is a completely solvable Lie group. This implies (see \cite{S}, Theorem 5):

\begin{theo} \label{isomequiv}
Two lattices $L$ and $L'$ of $\Osc$ are isomorphic if and only if their underlying discrete groups are isomorphic.
\end{theo}

Due to this result we will only use the term `isomorphic' and do not distinguish between `isomorphic as lattices' and `isomorphic as discrete groups'. Nevertheless, it is often useful to consider both kinds of isomorphisms since, depending on the situation, sometimes it is easier to consider isomorphisms of discrete groups and sometimes it is easier to consider isomorphisms of Lie groups.  

\section{Discrete split oscillator groups}

There are different definitions of what a discrete Heisenberg group should be. Here we will use the following one. It is chosen in such a way that every lattice in the (continuous) three-dimensional Heisenberg group is isomorphic to a discrete Heisenberg group as an abstract discrete group.
 A discrete Heisenberg group is a discrete group $H_{\rm d}$ isomorphic to 
\begin{equation}\label{dHeis}
H^r(\ZZ):=\la\, \alpha, \beta, \gamma \mid \ \alpha \beta \alpha^{-1} \beta^{-1} =\gamma^r,\ \alpha \gamma=\gamma \alpha,\ \beta\gamma=\gamma\beta\,\ra,
\end{equation}
for some $r\in\NN_{>0}$. For different $r$, these groups are non-isomorphic. 
\begin{de}
A discrete generalised oscillator group is a semi-direct product of a discrete Heisenberg group $H_{\rm d}$ with an infinite cyclic group generated by an element $\delta$ with respect to a homomorphism $\phi:\ZZ\cdot\delta\rightarrow \Aut(H_{\rm d})$ for which the map $\overline{\phi(\delta)}$ induced by $\phi(\delta)$ on $H_{\rm d} / Z(H_{\rm d} )\cong \ZZ^2$ is conjugate in $\GL(2, \RR)$ to a rotation or to a hyperbolic rotation.
\end{de}
If $\overline{\phi(\delta)}$ is conjugate to a rotation, i.e., to an element of $\SO(2)$, the group is an ordinary discrete oscillator group in the sense of \cite{FK}. Here we are interested in the case, where $\overline{\phi(\delta)}$  is conjugate to a hyperbolic rotation, i.e., to an element of $\SO^+(1,1)\setminus \{I_2\}$. This condition is satisfied if and only if $\overline{\phi(\delta)}$ is an element of 
$$\bB:=\{B\in\SL(2,\ZZ)\mid \tr B>2\}.$$ We will call the corresponding generalised oscillator groups \emph{split discrete oscillator groups}.

For $B=$\begin{small} $\begin{pmatrix}a&b\\c&d \end{pmatrix}$\end{small} $\in\bB$ and $l\in(\ZZ_r)^2$, we define
$$\dOsc(B,l):=H^r(\ZZ)\rtimes_{\phi} \ZZ\cdot\delta,$$
where $\phi(\delta)(\alpha)=\alpha^a\beta^c\gamma^{-l_2}$, $\phi(\delta)(\beta)=\alpha^b\beta^d\gamma^{l_1}$ and $\phi(\delta)(\gamma)=\gamma$.

\begin{lm}\label{L52}
The definition of $\dOsc(B,l)$ is correct, i.e., it does not depend on the representative of the congruence class of $l$ in $(\ZZ_r)^2$. 
\end{lm} 
\proof We show that there is an isomorphism $F:\dOsc(B, (l_1,l_2))\to\dOsc(B,(l_1+r,l_2))$ that is the identity on $H^r(\ZZ)$ and maps $\delta$ to $\delta\alpha^{-1}$. As above, let $\phi(\delta)$ be the automorphism of $H^r(\ZZ)$ that is determined by $B$ and $(l_1,l_2)$. We will shorten the notation and write just $\phi$ instead of $\phi(\delta)$. Analogously, let $\phi'$ be the automorphism determined by $B$ and $(l_1+r,l_2)$. We defined $F$ on generators and have to prove that $F$ maps relations to the identity. This is clear for relations within $H^r(\ZZ)$. Furthermore, we have  
$$F(\delta\alpha\delta^{-1}(\phi(\alpha))^{-1})=\delta \alpha^{-1}\alpha (\delta\alpha^{-1})^{-1}(\phi(\alpha))^{-1}=\delta\alpha \delta^{-1}(\phi'(\alpha))^{-1}=1
$$
since $\phi(\alpha)=\phi'(\alpha)$. Similarly, $\phi(\beta)=\phi'(\beta)\gamma^r$ implies
\begin{eqnarray*}
F(\delta\beta\delta^{-1}(\phi(\beta))^{-1})&=&\delta \alpha^{-1}\beta (\delta\alpha^{-1})^{-1}(\phi(\beta))^{-1}\\
&=&\delta\beta\delta^{-1}\gamma^{-r}(\phi(\beta))^{-1}\ =\ \delta\beta\delta^{-1}(\phi'(\beta))^{-1}=1.
\end{eqnarray*}
Analogously, there is an isomorphism $\dOsc(B, (l_1,l_2))\to\dOsc(B,(l_1,l_2+r))$ that is the identity on $H^r(\ZZ)$ and maps $\delta$ to $\delta\beta^{-1}$.
\qed

By definition, every discrete split oscillator group is isomorphic to a group $\dOsc(B,l)$ for some $B\in \bB$ and $l\in(\ZZ_r)^2$.  For different $(B,l)$, these groups are not necessarily non-isomorphic. In the next step, we clarify, what the isomorphism classes are. More exactly, we will parametrise the set $\cG$ of isomorphism classes of discrete split oscillator groups. 

Obviously, we have
$$\bB=\bigcup_{T=3}^\infty \bB_T, \quad \bB_T:=\{B\in \bB\mid \tr(B)=T\}.$$
We define an action of the group $\ZZ_2\times(\GL(2,\ZZ)\ltimes \ZZ^2)$ on $\bB_T\times (\ZZ_r)^2$ by  
\begin{eqnarray}
(K,m)\cdot(B,l)&:=& \big(KBK^{-1},K(l-(I_2-B^{-1})(m)) \big) \label{ga1}\\
\kappa \cdot(B,l)&:=& (B^{-1},-Bl) \label{ga2}
\end{eqnarray}
for $(K,m)\in \GL(2,\ZZ)\ltimes \ZZ^2$ and $\kappa=-1\in\ZZ_2$.
\begin{pr}\label{pre}
\begin{enumerate}
\item For $\dOsc(B,l)$, the numbers $r$ and $T:=\tr(B)$ depend only on the isomorphism class of $\dOsc(B,l)$.
\item Let $\cG_{r,T}$ be the set of isomorphism classes of discrete split oscillator groups with the same fixed $r$ and $T$. Then 
$$\cG=\coprod_{r=1}^\infty \  \coprod_{T=3}^\infty \cG_{r,T}$$
and the map $\bB_T\times (\ZZ_r)^2\ni(B,l) \mapsto \dOsc(B,l)$ induces a bijection 
\begin{equation}\label{bij}
\begin{array}{ccccc}  &&\hspace{-1em}\bB_T\times (\ZZ_r)^2&&\\[-2ex] &\hspace{-1em}\big\backslash &&\longrightarrow&\cG_{r,T}.\\[-2ex] \ZZ_2\times(\GL(2,\ZZ)\ltimes \ZZ^2)&&&&\end{array}
\end{equation}
\end{enumerate}
\end{pr}
\proof The commutator subgroup $\Gamma_0:=[\Gamma,\Gamma]$ of $\Gamma:=\dOsc(B,l)$ is the discrete Heisenberg group $H^r(\ZZ)$, hence $r$ is uniquely determined by the group structure of $\dOsc(B,l)$, thus it is invariant under isomorphisms. The same holds for $T$. Indeed, choose an element $\Delta$ of $\Gamma$ such that $\Gamma$ is generated by $\Gamma_0$ and $\Delta$. Then the action of $\Delta$ on $\Gamma_0/Z(\Gamma_0)\cong \ZZ^2$ equals that of $\delta$ or of $\delta^{-1}$. Thus the trace of this action equals $T$ and therefore $T$ is determined by the group structure of $\Gamma$. This proves the first assertion.   

For $K=$\begin{small} $\begin{pmatrix}k_1&k_2\\k_3&k_4 \end{pmatrix}$\end{small} $\in\GL(2,\ZZ)$, $m=(m_1,m_2)\in\ZZ^2$ and $\kappa\in\ZZ_2$, we define isomorphisms 
\begin{eqnarray*}
&F_{K}:& \ \dOsc(B,l) \longrightarrow \dOsc(\bar B,\bar l), \ (\bar B,\bar l)=(K,0)\cdot(B,l),\\
&F_{m}:& \ \dOsc(B,l) \longrightarrow \dOsc(\bar B,\bar l), \ (\bar B,\bar l)=(I,m)\cdot(B,l),\\
&F_{\kappa}:& \ \dOsc(B,l) \longrightarrow \dOsc(\bar B,\bar l), \ (\bar B,\bar l)=\kappa\cdot(B,l),
\end{eqnarray*}
by
\begin{eqnarray*}
&F_{K}:& \alpha \longmapsto\alpha^{k_1}\beta^{k_3},\quad \beta\longmapsto \alpha^{k_2}\beta^{k_4},\quad \gamma \longmapsto \gamma^{\det(K)} ,\quad \delta\longmapsto \delta\\
&F_{m}:&   \alpha\longmapsto \alpha\gamma^{-m_2},\quad \beta\longmapsto \beta\gamma^{m_1},\quad \gamma \longmapsto\gamma ,\quad \delta\longmapsto\delta \\
&F_{\kappa}:& \alpha\longmapsto \alpha,\quad \beta\longmapsto \beta,\quad \gamma \longmapsto\gamma ,\quad \delta\longmapsto \delta^{\kappa}.
\end{eqnarray*}
Then $F_\kappa\circ F_K\circ F_m: \dOsc(B,l) \rightarrow \dOsc(\bar B,\bar l)$,  $(\bar B,\bar l)=\kappa\cdot (K,m)\cdot(B,l)$, is an isomorphism. In particular, the map (\ref{bij}) is well-defined. Obviously, it is surjective. Let us prove that it is also injective. Assume that $\dOsc(B,l)$ and $\dOsc(\bar B,\bar l)$ are isomorphic and let $F:\dOsc(B,l)\rightarrow\dOsc(\bar B,\bar l)$ be an isomorphism. Since $H^r(\ZZ)$ is the commutator group of $\dOsc(B,l)$ and of $\dOsc(\bar B,\bar l)$, the isomorphism $F$ maps $H^r(\ZZ)\subset \dOsc(B,l)$ to $H^r(\ZZ)\subset \dOsc(\bar B,\bar l)$. Thus $F$ restricted to $H^r(\ZZ)$ is an automorphism of the discrete Heisenberg group. Consequently, $F$ is of the form
$$F:\ \alpha \longmapsto\alpha^{k_1}\beta^{k_3}\gamma^{-m_2\det(K)},\ \beta\longmapsto \alpha^{k_2}\beta^{k_4}\gamma^{m_1\det(K)},\ \gamma \longmapsto \gamma^{\det(K)} ,\ \delta\longmapsto (\delta\cdot h)^\kappa,$$
where $K\in\GL(2,\ZZ)$ is defined as above, $m=(m_1,m_2)\in\ZZ_2$, $h\in H^r(\ZZ)$ and $\kappa=\pm1$. Hence $F=F_m\circ F_K\circ F_\kappa\circ F_0$, where $F_0$ is the identity on $H^r(\ZZ)$ and maps $\delta$ to $\delta\cdot h$. Then $F_0$ maps $\dOsc(B,l)$ to $\dOsc(B,l')$ for some $l'$ with $l'\equiv l$ mod $r$.   We have seen this in the proof of Lemma~\ref{L52} for $h=\alpha^{-1}$ and $h=\beta^{-1}$, hence this is true for arbitrary $h\in H^r(\ZZ)$.  Consequently, 
$$F:\ \dOsc(B,l)  \overset{F_0}{\longrightarrow} \dOsc(B,l')\overset{F_m\circ F_K \circ F_\kappa}{\longrightarrow}\dOsc(\bar B,\bar l),$$
which gives $(\bar B,\bar l)=(I,m)\cdot(K,0)\cdot\kappa\cdot (B,l')$. Thus $(\bar B,\bar l)$ and $(B,l')=(B,l)\in\bB_T\times (\ZZ_r)^2$ are in the same orbit of $\ZZ_2\times(\GL(2,\ZZ)\ltimes \ZZ^2$.
\qedohne
\begin{de}For $B\in \SL(2,\Z)$, we say that $M \in \GL(2, \ZZ)$ is a symmetry of $B$ if $MBM^{-1}=B$ and we say that $M\in\GL(2,\ZZ)$ is a reversing symmetry if $MBM^{-1} =B^{-1}$. If $B$ admits a reversing symmetry, it is called reversible. Let $\cS(B)\subset \GL(2,\ZZ)$ denote the group of all symmetries of $B$ and $\cR(B)\subset \GL(2,\ZZ)$ the group of all symmetries and reversing symmetries of $B$. 
\end{de}
By definition, $\cR(B)$ is equal to $\cS(B)$ if $B$ is not reversible. Otherwise, $\cR(B)$ is generated by $\cS(B)$ and a single reversing symmetry of $B$.
\begin{re}\label{gensymm}
Take $B\in\bB_T$ and let $\lambda>1$ and $\lambda^{-1}$ be the (real) eigenvalues of $B$, i.e., $\lambda+\lambda^{-1}=T$. We denote by $K$ the field $\QQ(\lambda)$ and by $\cO_K$ the ring of integers in $K$. Then $\cS(B)$ is a subgroup of  $\{xI_2 + yB\mid x,y\in\QQ,\ x+y\lambda\in \cO_K^*\}$. The latter group is isomorphic to $\ZZ_2\times \ZZ$, 
where the generators correspond to $-1$ and to the fundamental unit $\eps$ in $\cO_K$. Hence $\cS(B)$ is also isomorphic to $\ZZ_2\times \ZZ$. Its generators correspond to $-1$ and to the smallest power $\eps^k=x_0+y_0\lambda$, $k\in\NN$, of the fundamental unit such that $x_0I_2+y_0B$ is an integer matrix. 
\end{re}
Fix $B\in\cB_T$. We define an action of $\cR(B)\ltimes \ZZ^2$ on $(\ZZ_r)^2$ by 
$$(K,m)\cdot l:=\left\{ \begin{array}{ll} K(l-(I_2-B^{-1})(m)),& \mbox{if $K\in \cS(B)$,}\\[1ex]  -B^{-1}K(l-(I_2-B^{-1})(m)),&\mbox{if $K\in\cR(B)\setminus \cS(B).$} \end{array}\right.$$ 
\begin{co} \label{cor} For $B\in\bB_T$, the groups
$\dOsc(B,l)$ and $\dOsc(B,l')$ are isomorphic if and only if $l$ and $l'$ are in the same  $\cR(B)\ltimes \ZZ^2$-orbit.
\end{co}
\proof By Prop.~\ref{pre}, $\dOsc(B,l)$ and $\dOsc(B,l')$ are isomorphic if and only if there exists an element $(\kappa,K,m)\in\ZZ_2\times(\GL(2,\ZZ)\ltimes \ZZ^2)$ such that $(\kappa,K,m)\cdot(B,l)=(B,l')$, where the group action is defined by (\ref{ga1}) and (\ref{ga2}). This implies that either $\kappa=1$ and $K$ is in the symmetry group of $B$ or $\kappa=-1$ and $K$ is a reversing symmetry of $B$. Now the assertion follows from (\ref{ga1}) and (\ref{ga2}). \qed

In Section~\ref{S7},  Corollary~\ref{Bl}, we will see that for each discrete oscillator group there is a lattice in $\Osc$ that is isomorphic to this group. Moreover, we will prove that every lattice in $\Osc$ is isomorphic to a discrete split oscillator group and we will determine this group.
\section{Parametrisation of the set of lattices}\label{S6}
\subsection{Arbitrary lattices in the split oscillator group}\label{S61}
The aim of this subsection is to parametrise the set of all lattices in $\Osc$. The main result is Proposition~\ref{Ppar}.

For our first observations on lattices in $\Osc$, we will use the following facts.  
\begin{itemize} 
\item[-] Let $G$ be a connected solvable Lie group  and $N$ be its maximal connected normal nilpotent Lie subgroup. If $\Gamma$ is a lattice in $G$, then $\Gamma \cap N$ is a lattice in $N$ (see \cite{Ra72}, Cor.\,3.5.).
\item[-] If $\Gamma_0$ is a lattice in the simply-connected nilpotent Lie group $N$, then the image of $\Gamma_0$ in the factor group $N/[N,N]$ is also a lattice (see {\rm \cite{Mal}}, Thm.\,4). 
\end{itemize} 
Now let $L$ be a lattice in $\Osc$. Then the first fact implies that $L\cap H$ is a lattice in the Heisenberg group. By the second fact, the image $\Lambda$ of $L\cap H$ in $H/Z(H)\cong\RR^2$ is a lattice in $\RR^2$. Let $p$ denote the projection of $\Osc=H\rtimes \RR$ onto the $\RR$-factor. Then the subgroup $p(L)$ of $\RR$ is discrete. Indeed, $e^{tA}(\Lambda)=\Lambda$ holds for each $t\in p(L)$. If $p(L)$ were dense in $\RR$, then $e^{tA}(\Lambda)=\Lambda$ would hold for all $t\in\RR$. This would imply that $e^{tA}$ is equal to the identity for all $t\in\RR$, which is a contradiction.

We want to use these observations to define a positive integer $r(L)$ and a positive real number $s(L)$ that is contained in the discrete set 
$$\SSS:=\{ s \in\RR_{>0}\mid e^{s}+e^{-s}\in\NN\}.$$ 
Since the intersection $L\cap H$ is a lattice in the Heisenberg group, it is isomorphic to some $\Gamma_r$ as defined in~(\ref{Gammar}). We put $r(L):=r$. Furthermore, we have seen that $p(L)\subset \RR$ is discrete, hence $p(L)=s\cdot\ZZ$ for a real number $s>0$. We define $s(L):=s$. Let $\delta$ be an element of $L$ such that $p(L)=s$. Then $\delta$ acts on $H/Z(H)\cong\RR^2$ with eigenvalues $e^{\pm s}$. Since this action preserves $\Lambda$, the trace $T$ of the matrix of this action is an integer, hence $T=e^s+e^{-s}\in\NN$.

We denote by $\cL$ the set of lattices in $\Osc$. Furthermore, for given $r\in\NN_{>0}$ and $s \in\SSS$, we define 
$$\cL_{r,s}:=\{ L\in \cL\mid r(L)=r,\ s(L)=s\}.$$ 

Next we turn to the set of parameters. We will use the set $\cB_r$, see Definition~\ref{df}. We put
\begin{eqnarray*}
{\bar\cP}_{r,s}&:=& \{(P,\xi)\mid P\in \GL(2,\RR),\, (\xi,B:=P^{-1}e^{sA}P)\in \cB_r\},\\
\cP_{r,s}&:=& \{(\eta, P,z,\xi)\mid \eta\in\RR^2,\,(P,\xi)\in \bar\cP_{r,s},\, z\in\RR\}, 
\end{eqnarray*}
and
$\displaystyle \cP=\coprod_{r=1}^\infty\coprod_{s\in{\Bbb S}}\cP_{r,s}.$
\begin{lm} Let $(\eta,P,z,\xi)$ be in $\cP_{r,s}$. Then the conjugation by $\delta:=(\det(P)z,P\xi,s)=(\bar F_{P}(z,\xi),s)$ stabilises the lattice $\bar F_{P}(\Gamma_r)\subset H$.
\end{lm}
\proof By definition of $\cP_{r,s}$, the pair $(\xi,B:=P^{-1}e^{sA}P)$ is in $\cB_r$. Hence, $(0,P\xi,s)$ stabilises $\bar F_{P}(\Gamma_r)$ as we have seen in Example~\ref{LxiP}. Thus this is also true for $\delta$ since $(\det(P)z,0,0)$ is in the centre of $\Osc$. \qed

The lemma shows that we can define lattices in $\Osc$ in the following way.
\begin{de}
Let $(\eta,P,z,\xi)$ be in $\cP_{r,s}$. We denote by $L_{r,s}(\eta,P,z,\xi)$ the image under $F_{\eta}$ of the lattice generated by $\bar F_{P}(\Gamma_r)\subset H$ and $(\bar F_{P}(z,\xi),s)$, i.e.,
$$ L_{r,s}(\eta, P,z,\xi):=F_\eta\big(\big\langle\, \bar F_{P}(\,\Gamma_r), (\bar F_{P}(z,\xi),s)\, \big\rangle \big).$$
\end{de}
Although we will not use explicit formulas for generators in this paper, we want to include them here. The lattice  $L_{r,s}(\eta, P,z,\xi)$ is generated by 
$$\textstyle(\frac1r\det P,0,0),\quad \big(\omega(\eta,Pe_1),Pe_1,0\big),\quad \big(\omega(\eta,Pe_2),Pe_2,0\big)$$
and
$$\textstyle\big( z\det P-\frac12\omega(P\xi+\eta,\eta+e^{sA}\eta),\, P\xi+\eta-e^{sA}\eta,\,s\big).$$
\begin{re}
For $r\in\NN_{>0}$, $s\in\SSS$ and $(P,\xi)\in {\bar{\cP}}_{r,s}$, the lattice $L_{r,s}(P,\xi)$ defined in Example~\ref{LxiP} equals 
$L_{r,s}(0,P,0,\xi).$
\end{re}
\begin{pr}
For each $r\in\NN_{>0}$ and $s \in\SSS$, the map 
\begin{equation} \label{sur}
\cP_{r,s} \longrightarrow \cL_{r,s}, \quad (\eta,P,z,\xi) \longmapsto L_{r,s}(\eta,P,z,\xi)
\end{equation}
is surjective. In particular, this defines a surjection from $\cP$ to $\cL$.
\end{pr}
\proof Let $L$ be a lattice in $\Osc$ such that $r(L)=r$ and $s(L)=s$. Then $L\cap H$ is a lattice in the Heisenberg group and isomorphic to $\Gamma_r$. Consequently, $L\cap H=\bar F_\eta\bar F_P (\Gamma_r)$ for some $P\in\GL(2,\RR)$ and $\eta\in\RR^2$. Now choose an element $\delta$ of $L$ such that $L$ is generated by $L\cap H$ and $\delta$. By assumption, the action of $\delta$ on $H/Z(H)\cong\ZZ_2$ is conjugate to $e^{sA}$. Hence we can write $\delta$ as $F_\eta\big((\bar F_{P}(z,\xi),s)\big)$ for suitable $z\in\RR$ and $\xi\in\RR^2$. It remains to prove that $(P,\xi)\in \bar\cP_{r,s}$. Since conjugation by $\delta$ stabilises $L\cap H=\bar F_\eta\bar F_P (\Gamma_r)$, conjugation by $(\bar F_{P}(z,\xi),s)$ stabilises $\bar F_P (\Gamma_r)$. Hence also conjugation by $(0,P\xi,s)$ stabilises $\bar F_P (\Gamma_r)$. Now the same calculation as in Example~\ref{LxiP} shows that $(\xi,B:=P^{-1}e^{sA}P)$ is in $\cB_r$.\qed

The map defined in (\ref{sur}) is not injective. Therefore, our next aim is to form a quotient of $\cP_{r,s}$ such that the induced map becomes bijective. By Definition~\ref{df}, $\cB_r$ acts on $\Gamma_r$. Thus we can consider the group $\Gamma_r\rtimes \cB_r$. We want to define an action of this group on the right of $\cP_{r,s}$. Let $(k,l,\zeta,K)$ be an element of $\Gamma_r\rtimes \cB_r$, where $(k,l)\in\Gamma_r\subset H$ and $(\zeta,K)\in \cB_r$. Then we put
\begin{equation}\label{ga}
(\eta, P, z, \xi)\cdot (k,l,\zeta, K)=(\tilde\eta, \tilde P, \tilde z, \tilde \xi),
\end{equation}
where 
\begin{eqnarray}
(\tilde \eta,\tilde P)& = & (\eta,P)\cdot (\zeta, K) \label{tilde1}\\
(\tilde z, \tilde \xi) &=& \bar F_{K^{-1}}\big( (0,-\zeta)(k,l)^{-1}(z,\xi)(0,B\zeta)\big), \label{tilde2}
\end{eqnarray}
where, as usual, $B=P^{-1}e^{sA}P$. The multiplication in (\ref{tilde1}) is the multiplication in the group $\RR^2\rtimes \GL(2,\RR)$ and the multiplication in (\ref{tilde2}) is the multiplication in $H$.
\begin{re}
Evaluating the second component in (\ref{tilde2}), we obtain the equation 
\begin{equation}\label{tildexi}
\tilde \xi= K^{-1}(\xi -l+(B-I)\zeta).
\end{equation}
\end{re}
\begin{lm} Equations (\ref{ga}) -- (\ref{tilde2}) define an action of $\Gamma_r\rtimes \cB_r$ on $\cP_{r,s}$.
\end{lm}
\proof We show that $(\tilde\eta, \tilde P, \tilde z, \tilde \xi)$ is again in $\cP_{r,s}$. We have to check that $(\tilde \xi,\tilde B)\in \cB_r$ holds for $\tilde \xi$ as given by (\ref{tildexi}) and $\tilde B=\tilde P^{-1}e^{sA}\tilde P=K^{-1}BK$. This follows from
$$(\tilde \xi,\tilde B)=(K^{-1}(\xi -l+(B-I)\zeta),K^{-1}BK)=(\zeta,K)^{-1}(-l,I)(\xi,B)(\zeta,K)$$
in $\RR^2\rtimes \GL(2,\RR)$ and the fact that each factor on the right hand side is in $\cB_r$.
It is now a straightforward calculation to check that (\ref{ga}) is indeed a group action.
\qedohne
\begin{pr}\label{Ppar} \begin{enumerate} 
\item The sets $\cL_{r,s}$ are pairwise disjoint, i.e., $\cL=\coprod_{r=1}^\infty\coprod_{s\in{\Bbb S}}\cL_{r,s}$. 
\item 
The map defined in (\ref{sur}) induces a bijection $\cP_{r,s}/ (\Gamma_r\rtimes \cB_r)\rightarrow \cL_{r,s}$.
\end{enumerate}
\end{pr}
\proof The first assertion is obvious since, for $L\in\cL$, the numbers $r(L)$ and $s(L)$ are determined by $L$. Suppose that $L_{r,s}(\eta,P,z,\xi)=L_{r,s}(\tilde \eta,\tilde P,\tilde z,\tilde \xi)$. Then $F_{\eta} \bar F_{P}(\,\Gamma_r)=F_{\tilde\eta} \bar F_{\tilde P}(\,\Gamma_r)$ holds, which is equivalent to 
$$\Gamma_r=\bar F_{P^{-1}}F_{-\eta}F_{\tilde \eta} \bar F_{\tilde P}(\,\Gamma_r)=F_{P^{-1}(\tilde \eta-\eta)} \bar F_{P^{-1} \tilde P}(\,\Gamma_r),$$
thus to $(\zeta,K)\in \cB_r$ for $\zeta:=P^{-1}(\tilde \eta-\eta)$ and $K:=P^{-1}\tilde P$. We obtain $(\tilde\eta,\tilde P)=(\eta,P)\cdot(\zeta,K)$.
Furthermore,
\begin{eqnarray*}
F_{\tilde\eta}(\bar F_{\tilde P}(\tilde z,\tilde\xi),s)&=&F_{\eta}\big((\bar F_{P}(k,l)^{-1},0)\cdot(\bar F_{P}(z,\xi),s)\big)\\
&=&\big(\bar F_{\eta}\bar F_{P}((k,l)^{-1}(z,\xi)),0\big)\cdot F_{\eta}(0,0,s)
\end{eqnarray*}
is equivalent to 
\begin{eqnarray*}
(\bar F_{\tilde P}(\tilde z,\tilde\xi),s)&=&\big(\bar F_{-P\zeta}\bar F_{P}((k,l)^{-1}(z,\xi)),0\big)\cdot F_{-P\zeta}(0,0,s)\\
&=&\big(\bar F_{-P\zeta}\bar F_{P}((k,l)^{-1}(z,\xi)),0\big)\cdot (0, -P\zeta,0)\cdot(0,0,s)\cdot(0,P\zeta,0)\\
&=&\big(\bar F_{-P\zeta}\bar F_{P}((k,l)^{-1}(z,\xi)),0\big)\cdot (0, -P\zeta,0)\cdot(0,e^{sA}P\zeta,0)\cdot(0,0,s)
\end{eqnarray*}
thus to the equality 
\begin{eqnarray*}
(\tilde z,\tilde\xi)
&=&\bar F_{\tilde P^{-1}} \big(\bar F_{-P\zeta}\bar F_{P}((k,l)^{-1}(z,\xi)) (0,-P\zeta)(0,e^{sA}P\zeta)\big)\\
&=&\bar F_{\tilde P^{-1}} \big(\bar F_{P}\bar F_{-\zeta}((k,l)^{-1}(z,\xi)) (0,-P\zeta)(0,e^{sA}P\zeta)\big)\\
&=&\bar F_{\tilde P^{-1}} \big(\bar F_{P}((0,-\zeta)(k,l)^{-1}(z,\xi)(0,P^{-1}e^{sA}P\zeta))\big)\\
&=&\bar F_{K^{-1}}\big((0,-\zeta)(k,l)^{-1}(z,\xi)(0,B\zeta)\big)
\end{eqnarray*}
in the Heisenberg group.
\qed
\subsection{Normalised lattices and their parametrisation}
For a lattice $L$ in $\Osc$, let $\bar L_0$ denote the projection of $L_0:=L\cap H$ to $H/Z(H)\cong \RR^2$.
\begin{de} 
The lattice $L$ is called normalised if $\bar L_0$ is a normalised lattice in $\RR^2$, i.e., if $\bar L_0\subset \RR^2$ has covolume one with respect to the standard metric of~$\RR^2$.
\end{de}
Let $\cL^0$ denote the set of normalised lattices in $\Osc$.
\begin{pr}
The map 
$$\cL^0\times \RR_{>0}\longrightarrow \cL,\quad (L,a)\longmapsto F_{S}(L),\ S=aI_2,$$
is a bijection.
\end{pr}
\proof We prove the assertion by determining the inverse map. Take $L\in\cL$. Let $a(L)$ be the square root of the covolume of $L_0\subset\RR^2$. Then the inverse map is given by $L\mapsto (F_S (L), a(L))$ for $S=a(L)^{-1}I_2$. \qed

We can also obtain a parametrisation of $\cL^0$. We put 
$$\cP_{r,s}^0:= \{(\eta, P,z,\xi)\in \cP_{r,s}\mid \det P=\pm1\},\quad \cP^0:=\coprod_{r=1}^\infty\coprod_{s\in{\Bbb S}}\cP_{r,s}^0\,.$$
In analogy to Prop.~\ref{Ppar}, we obtain:
\begin{pr}\label{Ppar0} 
The set $\cL^0$  equals the disjoint union $\coprod_{r=1}^\infty\coprod_{s\in{\Bbb S}}\cL^0_{r,s}$, where $\cL^0_{r,s}$
is in bijection with $\cP^0_{r,s}/ (\Gamma_r\rtimes \cB_r)$.
\end{pr}
\section{Isomorphism classes of lattices}\label{S7}

In Section~\ref{S6}, we defined a lattice $L_{r,s}(\eta,P,z,\xi)$ for every $(r,s)\in \NN_{>0}\times\SSS$ and every $(\eta,P,z,\xi) \in\cP_{r,s}$. This gave us a parametrisation of the set $\cL$ of all lattices, see Prop.~\ref{Ppar}. Now we want to decide for which parameters $(r,s,\eta,P,z,\xi)$ we obtain isomorphic lattices. In the end of this section we show that each lattice $L_{r,s}(\eta,P,z,\xi)$ is isomorphic to a discrete split oscillator group and we determine this group. 

Recall that $L_{r,s}(P,\xi)=L_{r,s}(0,P,0,\xi)$.
\begin{pr} \label{Piso} \begin{enumerate}
\item The lattice $L_{r,s}(\eta,P,z,\xi)$ is isomorphic to $L_{r,s}(P,\xi)$. 
\item
The lattices $L_{r,s}(P,\xi)$ and $L_{r',s'}(P',\xi')$ are isomorphic if and only if $(r,s)=(r',s')$ and if there exists an element $(\zeta,K)\in\cB_r$ such that 
\begin{enumerate}
\item $B'=KBK^{-1}$, and
\item $\xi'-K\xi+(B'-I)\zeta\in \ZZ^2$,
\end{enumerate}
or \begin{enumerate}
\item $B'=KB^{-1}K^{-1}$, and
\item $\xi' +KB^{-1}  \xi+(B'-I)\zeta\in \ZZ^2$. 
\end{enumerate}
\end{enumerate}
\end{pr}
\proof Obviously, $F_{-\eta}(L_{r,s}(\eta,P,z,\xi)=L_{r,s}(0,P,z,\xi)$. Moreover, put $u:=-\det(P)z/s$. Then
$$F_u(L_{r,s}(0,P,z,\xi))=F_u\big(\big\langle \bar F_P(\Gamma_r), (\bar F_P(z,\xi),s)\big\rangle\big)=\big\langle \bar F_P(\Gamma_r), F_u(\bar F_P(z,\xi),s)\big\rangle$$
and 
$$F_u(\bar F_P(z,\xi),s)=F_u(\det(P)z,P\xi,s)=(\det(P)z+us,P\xi,s)=(0,P\xi,0)=(\bar F_P(0,\xi),s),$$
which implies $F_u(L_{r,s}(0,P,z,\xi))=L_{r,s}(0,P,0,\xi)$. This proves the first assertion.

For any lattice $L\in\cL$, the numbers $s(L)$ and $r(L)$ only depend on the isomorphism class of $L$, which gives the first assertion in Item (ii). It remains to study under which conditions $L_{r,s}(P,\xi)$ and $L_{r,s}(P',\xi')$ are isomorphic. An arbitrary isomorphism of $\Osc$ is of the form $F:=F_\eta\circ F_S\circ F_u$, see Subsection~\ref{S22}. We have 
\begin{eqnarray*}
F(L_{r,s}(P,\xi))&=&(F_\eta\circ F_S\circ F_u)\big(\big\langle \bar F_P(\Gamma_r), (\bar F_P(0,\xi),s)\big\rangle\big)\\
&=&(F_\eta\circ F_S)\big(\big\langle \bar F_{P}(\Gamma_r), (us,P\xi, s)\big\rangle\big)\\
&=&F_\eta\big(\big\langle \bar F_{SP}(\Gamma_r), (us\det(S),SP\xi, \mu s)\big\rangle\big),
\end{eqnarray*}
where $\mu$ is the sign of $\det S$.

If $\mu=1$, then $$(us\det(S),SP\xi, \mu s)=(\bar F_{SP}(us\det(P)^{-1},\xi),s).$$
If $\mu=-1$, then
\begin{eqnarray*}
(us\det(S),SP\xi, \mu s)&=&(us\det(S),SP\xi, -s)\\
&=&(-us\det(S),-e^{sA}SP\xi,s)^{-1}\\
&=&\big(\bar F_{SP}(-us\det(P)^{-1},-(SP)^{-1}e^{sA}SP\xi),s\big)^{-1}\\
&=&\big(\bar F_{SP}(-us\det(P)^{-1},-P^{-1}e^{-sA}P\xi),s\big)^{-1}\\
&=& \big(\bar F_{SP}(-us\det(P)^{-1},-B^{-1}\xi),s\big)^{-1}
\end{eqnarray*}
Thus
$$F(L_{r,s}(P,\xi))=\left\{\begin{array}{ll} L_{r,s}(\eta,P,us\det(P)^{-1},\xi)\,, &\mbox{if } \mu=1,\\[1ex]
 L_{r,s}(\eta,SP,-us\det(P)^{-1},-B^{-1}\xi  )\,, &\mbox{if } \mu=-1.
\end{array}\right.$$
Consequently, $L_{r,s}(P,\xi)$ and $L_{r,s}(P',\xi')$ are isomorphic if and only if there are $u,\eta,S,k,l,K,\zeta$ such that 
$$(0,P',0,\xi')\cdot (k,l,\zeta, K) =\left\{\begin{array}{ll}(\eta,SP,us\det(P)^{-1},\xi)\,,&\mbox{if } \mu=1,\\[1ex] (\eta,SP,-us\det(P)^{-1}, -B^{-1}\xi )\,, &\mbox{if } \mu=-1.\end{array}\right.$$
This is equivalent to the existence of $S,l,\zeta,K$ such that $P'K=SP$ and 
$$K^{-1}(\xi'-l+(B'-I)\zeta)= \left\{\begin{array}{ll}\xi\,, &\mbox{if } \mu=1,\\[1ex]
 -B^{-1}\xi \,, &\mbox{if } \mu=-1,\end{array}\right.$$
 see (\ref{tildexi}).
The existence of $S$ is equivalent to $P'KP^{-1}\in\cS$, that is, to the condition $P'KP^{-1}e^{sA}PK^{-1}P'^{-1}= e^{\mu sA} $. The last equation can be rewritten as $KBK^{-1}=(B')^{\mu}$. Then the existence of $l$ is equivalent to Condition (b) in the Proposition. 
\qedohne
\begin{co}\label{Bl}
The lattice $L_{r,s}(\eta,P,z,\xi)$ is isomorphic to the discrete group $\dOsc(B,l)$ for 
\begin{equation} 
B=P^{-1}e^{sA}P=:\begin{pmatrix}a&b\\c&d \end{pmatrix}, \quad l=rB^{-1}\xi+{\textstyle \frac r2} {-bd\choose ac}. \label{Pass}
\end{equation}
In particular, every lattice in $\Osc$ is isomorphic to a discrete split oscillator group. 

Moreover, for every discrete oscillator group there is a lattice in $\Osc$ that is isomorphic to this group. More exactly, given a discrete split oscillator group $\dOsc(B,l)$, choose $P$ such that $B=P^{-1}e^{sA}P$ and put 
$$
\xi:={\textstyle\frac1r}Bl -{\textstyle \frac 12} B{-bd\choose ac}, \mbox{ where } B=\begin{pmatrix}a&b\\c&d \end{pmatrix}.$$
Then the lattice $L_{r,s}(P,\xi)$ is isomorphic to $\dOsc(B,l)$.
\end{co}
\proof By Prop.~\ref{Piso}, we know that $L_{r,s}(\eta,P,z,\xi)$ is isomorphic to the lattice $L_{r,s}(P,\xi)$, which is generated by  $\alpha:=(0,Pe_1,0)$, $\beta:=(0,Pe_2,0)$, $\gamma:=(\frac1r\det P,0,0)$ and $\delta:=(0,P\xi,s)$. The elements $\alpha, \beta$ and $\gamma$ generate a discrete group isomorphic to $H^r(\ZZ)$. Then
\begin{eqnarray*}
\delta \alpha\delta^{-1}&= &F_{P\xi}F_{e^{sA}}(0,Pe_1,0)\ =\ \big(\omega(P\xi, e^{sA}Pe_1),e^{sA}Pe_1,0\big)\\
&=&\big(\det(P)\cdot\omega(\xi, Be_1),PBe_1,0\big).
\end{eqnarray*}
Furthermore,
\begin{eqnarray*}
\big(\det(P)\cdot\omega(\xi, Be_1),PBe_1\big)&=&\bar F_P\big(\omega(\xi, Be_1),Be_1\big)=\bar F_P\big(\xi_1 c-\xi_2 a, (a,c)^\top\big)\\
&=&\bar F_P\big( (\xi_1c-\xi_2a-\textstyle \frac12 ac,0)(0,e_1)^a(0,e_2)^c\big)\\
&=& (\textstyle\frac1r \det P,0)^{-l_2} (0,Pe_1)^a(0,Pe_2)^c
\end{eqnarray*}
for $l_2=-r(\xi_1c-\xi_2a)+\textstyle \frac r2 ac$. Together with a similar calculation for $\delta\beta\delta^{-1}$ this proves  Equation~(\ref{Pass}).

It remains to prove the last assertion. Note first that $(B,\xi)$ is in $\cB_r$. Indeed, if we define $\xi$ as above, then $\xi$ satisfies the second equation in (\ref{Pass}), which implies that (\ref{odd}) holds, which is equivalent to  $(B,\xi)\in\cB_r$, see the proof of Lemma~\ref{LBr}. Hence $L_{r,s}(P,\xi)$ is a lattice in $\Osc$. Now (\ref{Pass}) shows that $L_{r,s}(P,\xi)$ is isomorphic to $\dOsc(B,l)$.
\qedohne
\section{Classification}
We know that two lattices in $\Osc$ are isomorphic if and only if the underlying discrete split oscillator groups are isomorphic, see Section~\ref{S4}. In Prop.~\ref{pre} we gave a description of the set $\cG$ of isomorphism classes of discrete split oscillator groups. In this section, we want to make this description more explicit.
\subsection{Conjugacy and reduced cycles}\label{S81}
In Prop.~\ref{pre}, we proved that each $\cG_{r,T}$ is in bijection with a quotient space given by the left hand side of (\ref{bij}). If we want to use this bijection in order to study $\cG_{r,T}$, we need a good description of the quotient $\GL(2,\ZZ)\backslash \bB_T$, where $\GL(2,\ZZ)$ acts on $\bB_T\subset\SL(2,\ZZ)$ by conjugation. This is a classical problem and there are several methods to describe these conjugacy classes. Here we want to apply the one developed by Aicardi~\cite{A}. In the appendix we give a short introduction to this method, where we slightly modify it in order to adapt it to our needs. It relies on the notion of cycles in a finite subset $H^0_{\rm red}$ of $\bB_T$, which we want to introduce in the following.

Let $T\in\NN$, $T>2$, be fixed.
\begin{de}
Let $A$ and $B$ be elements of\ \,$\SL(n,\ZZ)$. We say that $A$ and $B$ are extendedly conjugate if $B$ is conjugate to $A$ or $A^{-1}$ in $\GL(n,\ZZ)$. The corresponding equivalence classes in $\SL(n,\ZZ)$ are called extended conjugacy classes.
\end{de}
We put 
$$\hat A:=\begin{pmatrix}1&1\\0&1 \end{pmatrix},\quad \hat B:=\begin{pmatrix}1&0\\1&1 \end{pmatrix}$$
and consider the finite set
$$H^0_{\rm red}=\left\{ \begin{pmatrix}a&b\\c&d\end{pmatrix} \in\bB_T \ \big|\  a,b,c,d>0,\ \max\{a,b,c,d\}\in\{a,d\} \right\}.$$
\begin{de}\label{redc}
A cycle of length $t$ in $H^0_{\rm red}$ is a cyclic sequence $[B_1,\dots, B_t]$ of $t>1$ distinct matrices in $H^0_{\rm red}$ such that $B_{i+1} = M_{i}B_{i}M_{i}^{-1}$ $(i = 1,\dots,t-1)$ and $B_1 = M_t B_t M_t ^{-1}$, where each of the matrices
$M_1,\dots,M_t$ is a power of $\hat A$ or $\hat B$. 
\end{de}
\begin{re} \label{expl} {\rm
\begin{enumerate}
\item  We consider the index $i$ of the elements in a cycle of length $t$ as an element of $\ZZ_t=\ZZ/t\ZZ$, i.e., we put $i+1=1$ if $i=t$ and $i-1=t$ if $i=1$.
\item We can describe cycles explicitly in the following way. Let
$B_i:=${\small $\begin{pmatrix} a_i&b_i\\ c_i&d_i \end{pmatrix}$}, $i\in\ZZ_t$, be the elements of a cycle in $H^0_{\rm red}$. Then
$a_i\not=d_i$. Moreover, $a_i>d_i$ holds if and only if $a_{i+1}<d_{i+1}$.

If $a_i<d_i$, then $M_i=\hat A^q$ for $$q=\left\{ \begin{array}{ll}
  \left[ \frac{d_i}{c_i}\right],    & \mbox{if } c_i\not=1,\\[1ex]
   d_i-1,  & else,
\end{array}\right.$$
and, if $a_i>d_i$, then $M_i=\hat B^q$ for $$q=\left\{ \begin{array}{ll}
  \left[ \frac{a_i}{b_i}\right],    & \mbox{if } b_i\not=1,\\[1ex]
   a_i-1,  & else.
\end{array}\right.$$
The statement follows from Corollary~\ref{Coapp} in the appendix and the remark that a reduced cycle in $H^0$ is the same as a cycle in $H^0_{\rm red}$.
\end{enumerate}
}\end{re}
The importance of cycles lies in the following 

{\bf Fact.} {\it The set $H^0_{\rm red}$ decomposes into pairwise disjoint cycles. Each ${\SL}(2,\ZZ)$-orbit in $\bB_T$ contains exactly one cycle in $H^0_{\rm red}$.}

For a proof, see Prop.~\ref{c} and Prop.~\ref{redunred}.
We introduce two further matrices
$$\hat U:=\begin{pmatrix} 0&1\\ 1&0 \end{pmatrix},\quad \hat W :=\begin{pmatrix} 0&-1\\ 1&0 \end{pmatrix}.$$
For $B\in H^0_{\rm red}$, we define $\bar B:=\hat U B \hat U^{-1}$.
If $\bz:=[B_1,\dots, B_t]$ is a cycle in $H^0_{\rm red}$, then $\bz^\top:=[B_1^\top,\dots, B_t^\top]$ and $\bar\bz:=[\bar B_1,\dots,\bar B_t]$ are also  cycles in $H^0_{\rm red}$. The cycles $\bz$ and $\bar \bz$ are in the same conjugacy class in $\GL(2,\ZZ)$. For $B\in\bz$, the inverse $B^{-1}$ and the elements of $\bz^\top$ are in the same conjugacy class in $\SL(2,\ZZ)$. Indeed, $\hat W B^{-1}\hat W^{-1}=B^\top$.
\begin{ex}\label{ExT} {\rm We want to decompose $H^0_{\rm red}$ into cycles for $T\in\{3,4,20\}$. Instead of writing the cycles in the formal way introduced above, we will display them as graphs. The elements of the cycle constitute the nodes of the graph and an arrow with label $\hat A^q$ or $\hat B^q$ stands for conjugation by the label.

For $T=3$, we have only one cycle:\\
{\small
\begin{tikzpicture}
\node (ll) at (1,0){$\bz=\bar\bz=\bz^\top:$};
\node (l) at (3,0) {$\begin{pmatrix} 2&1\\ 1&1 \end{pmatrix}$\ };
\node (r) at (5,0) {\ $\begin{pmatrix} 1&1\\ 1&2 \end{pmatrix}$};
\draw[->,bend left=45 ] (l) edge node [above] {$\hat B$} (r); 
\draw[<-,bend left=-45] (l) edge node [below] {$\hat A$} (r);
\end{tikzpicture}
}

For $T=4$, we have two cycles:\\
{\small
\begin{tikzpicture}
\node (ll) at (1,0){$\bz_1:\ $};
\node (l) at (2,0) {$\begin{pmatrix} 3&2\\ 1&1 \end{pmatrix}$\ };
\node (r) at (4,0) {\ $\begin{pmatrix} 1&2\\ 1&3 \end{pmatrix}$\ ,};
\draw[->,bend left=45 ] (l) edge node [above] {$\hat B$} (r); 
\draw[<-,bend left=-45] (l) edge node [below] {$\hat A^2$} (r);
\node (llb) at (7,0){$\bz_2=\bar \bz_1=\bz_1^\top:\ $};
\node (lb) at (9,0) {$\begin{pmatrix} 1&1\\ 2&3 \end{pmatrix}$\ };
\node (rb) at (11,0) {\ $\begin{pmatrix} 3&1\\ 2&1 \end{pmatrix}$\ .};
\draw[->,bend left=45 ] (lb) edge node [above] {$\hat A$} (rb); 
\draw[<-,bend left=-45] (lb) edge node [below] {$\hat B^2$} (rb);
\end{tikzpicture}
}

For $T=20$, there are six cycles:\\
{\small
\begin{tikzpicture}
\node (ll) at (1,0){$\bz_1:\ $};
\node (l) at (2.2,0) {$\begin{pmatrix} 19&18\\ 1&1 \end{pmatrix}$\ };
\node (r) at (4.5,0) {\ $\begin{pmatrix} 1&18\\ 1&19 \end{pmatrix}$\ ,};
\draw[->,bend left=45 ] (l) edge node [above] {$\hat B$} (r); 
\draw[<-,bend left=-45] (l) edge node [below] {$\hat A^{18}$} (r);
\node (llb) at (7,0){$\bz_2=\bar \bz_1=\bz_1^\top, $};
\end{tikzpicture}
}

{\small
\begin{tikzpicture}
\node (ll) at (1,0.9){$\bz_3:\ $};
\node (l) at (2.9,2) {$\begin{pmatrix} 18&7\\ 5&2 \end{pmatrix}$\ };
\node (r) at (5.2,2) {\ $\begin{pmatrix} 4&7\\ 9&16 \end{pmatrix}$\ };
\draw[->,bend left=45 ] (l) edge node [above] {$\hat B^2$} (r); 
\node (llb) at (9,0.9){$,\quad\ \bz_4=\bar \bz_2,\ \bz_5=\bz_2^\top,\ \bz_6=\bar\bz_2^\top .$};
\node (lu) at (3,0) {$\begin{pmatrix} 3&10\\ 5&17 \end{pmatrix}$\ };
\node (ru) at (5.2,0) {\ $\begin{pmatrix} 13&10\\ 9&7 \end{pmatrix}$\ };
\draw[<-,bend left=-45] (lu) edge node [below] {$\hat B$} (ru);
\draw[<-,bend left=-45] (l) edge node [left] {$\hat A^3$} (lu);
\draw[->,bend left=45] (r) edge node [right] {$\hat A$ } (ru);
\end{tikzpicture}
}

We have chosen $T=20$ as an example since this is the smallest number for which there exists a cycle $\bz$ such that $\bz$, $\bar \bz$, $\bz^\top$ and $\bar \bz^\top$ are pairwise different.
}\end{ex}
\subsection{Classification strategy} \label{S82}
In this subsection we give a method to determine the isomorphism classes of lattices in $\Osc$ or, equivalently, of all discrete split oscillator groups. 

Let us fix integers $T>2$ and $r>0$. For a fixed map $B\in\bB_T$, we define $\sim_B$ on $\ZZ_r^2$ as the smallest equivalence relation that contains the relation defined as follows. Two elements $l,l'\in\ZZ_r^2$ are related if  
\begin{itemize}
\item[(i)] there exists an element $m \in (\ZZ_r)^2$ such that $l'= l + (B-I_2)(m)$ or
\item[(ii)] there exists a map $K \in \mathcal{R}(B)$ such that $l' = Kl$.
\end{itemize}
Next we describe a procedure how we can find a complete set of representatives of $\cG_{r,T}$. The final result for $\cG$ will be formulated in Theorem~\ref{theo}.
\begin{enumerate} 
\item\textbf{Decompose $H^0_{\rm red}$ into cycles. } To find this decomposition, start with an element of $H^0_{\rm red}$ and determine its cycle by conjugating by powers of $\hat A$ and $\hat B$ according to Remark~\ref{expl}, 2. Then take one of the remaining elements of $H^0_{\rm red}$ and determine its cycle. Proceed until all cycles are determined. This happens after finitely many steps since $H^0_{\rm red}$ is finite. Instead of separate cycles now consider pairs of cycles $\bz,\bar\bz$, which may coincide. Choose a representative of each set $\bz\cup\bar\bz$ in this decomposition. Let us denote these representatives by $B_1,\dots,B_n$. 
\item \textbf{Determine the reversing symmetry group of each representative.}
\begin{enumerate}
\item\textbf{Determine the symmetry group of each representative. } For each representative $B_j$, $j=1,\dots,n$, determine an element $G_j \in \GL(2,\ZZ)$ such that $\cS(B_j) = \{ \pm G_j^q \mid q \in \ZZ\}$, see Remark~\ref{gensymm} for how to find $G_j$. If $B_j$ is not reversible, then $\cR(B_j)=\cS(B_j)$.
\item\textbf{Determine a reversing symmetry of each reversible representative. } An element $B$ of $H^0_{\rm red}$ is reversible, if and only if its cycle $\bz$ satisfies $\bz^\top\in\{\bz, \bar\bz\}$. For each reversible representative $B_j$ we find a reversing symmetry in the following way. If $B_j^\top\in\bz$, then we can read a matrix $M$ from the cycle $\bz$ such that $B_j^\top=MB_jM^{-1}$. Then $R_j:=\hat W^{-1}M$ is a reversing symmetry of $B_j$. If $B_j^\top$ is in $\bar \bz$, then $B_j^\top=M\bar B_jM^{-1}$, where $M$ can be read from $\bar \bz$. Now $R_j:=\hat W^{-1}M\hat U$ is a reversing symmetry. In both cases, $\cR(B_j)$ is generated by $\cS(B_j)$ and $R_j$. 
\end{enumerate}
\item\textbf{Determine the equivalence classes of $(\ZZ_r)^2$ for $\sim_{B_j}$, for each  $B_j$, $j=1,\dots,n$. }  
\end{enumerate}
Since $n$ and $B_1,\dots,B_n$ depend on $T$, we will use the notation $n(T)$ and $B_j^T$, $j=1,\dots,n(T)$ in the following.

\begin{theo} \label{theo} The map
\begin{eqnarray*}\coprod_{r=1}^\infty \  \coprod_{T=3}^\infty\ \coprod_{j=1}^{n(T)}\ (\ZZ_r)^2/_{\sim_{B_j^T}} \longrightarrow \cG,\quad
[l]\longmapsto  \dOsc(B_j^T,l)
\end{eqnarray*}
is a bijection onto the set $\cG$ of isomorphism classes of discrete split oscillator groups.
\end{theo}
\proof The statement almost follows from Proposition~\ref{pre} and Corollary~\ref{cor}. It remains to show that the quotient $(\cR(B)\ltimes \ZZ^2)\backslash (\ZZ_r)^2$ is in bijection with the set of equivalence classes $(\ZZ_r)^2/_{\sim_{B}}$ for each $B:=B_j^T$. Suppose first that $l$ and $l'$ are in the same $(\cR(B)\ltimes \ZZ^2)$-orbit. Then we have $l'=K(l-(I_2-B^{-1})(m))$ for some $K\in\cS(B)$ and $m\in\ZZ^2$ or $l'=-B^{-1}K(l-(I_2-B^{-1})(m))$ for some $K\in\cR(B)\setminus \cS(B)$ and $m\in\ZZ^2$. In both cases we obtain $l'\sim_B l-(I_2-B^{-1})(m)$ since $K$ and $-B^{-1}K$ are in $\cR(B)$. Furthermore, $l-(I_2-B^{-1})(m)=l-(B-I_2)(B^{-1}m)\sim_B l$, which proves $l'\sim_B l$. Conversely, suppose that $l'\sim_B l$. We may restrict to the relations defined by (i) and (ii) since these relations generate $\sim_B$. Thus we may assume that $l'= l + (B-I_2)(m)$ for some $m\in (\ZZ_r)^2$ or that $l'=Kl$ for some $K\in\cR(B)$. In the first case, we obtain $l'= l - (I_2-B^{-1})(-Bm)$. Thus, $l'=(I_2,-Bm)\cdot l$, which shows that $l$ and $l'$ are in the same orbit. In the second case, we have $l'=(K,0)\cdot l$ if $K\in\cS(B)$. If $K\in\cR(B)\setminus \cS(B)$, then also $-BK\in\cR(B)\setminus \cS(B)$ and we get $l'=(-BK,0)\cdot l$.
\qed

For given $T\in\NN$, $T\ge3$, we define $s_T>0$ by $e^{s_T}+e^{-s_T}=T$. For each $B_j^T$, we choose a matrix $P_j^T$ such that $B_j^T=(P_j^T)^{-1}e^{s_TA}P_j^T$. Corollary~\ref{Bl} now implies
\begin{co}
The map $$\coprod_{r=1}^\infty \  \coprod_{T=3}^\infty\ \coprod_{j=1}^{n(T)}\ \ZZ^2_r/_{\sim_{B_j^T}} \ni
[l]\longmapsto L_{r,s_T}(P_j^T,\xi),$$
where
$$ \xi:={\textstyle\frac1rB_j^T l - \frac 12} B_j^T{-bd\choose ac}\quad  \mbox{for}\quad  B_j^T = \begin{pmatrix}a&b\\c&d \end{pmatrix},$$
is a bijection onto the set of isomorphism classes of lattices in $\Osc$.
\end{co}

\subsection{Explicit results for small $T$}\label{S83}

For each discrete oscillator group $\dOsc(B,l)$, the number $T=\tr B\in\NN$ depends only on the isomorphism class of the group.
Here we want to classify all discrete oscillator groups with $T\le 7$ up to isomorphism. We proceed as described in Subsection~\ref{S82}. \\

\fbox{$T=3$}\quad We already know that there is only one cycle in $H^0_{\rm red}$, see Example~\ref{ExT}. It is represented by $B=$ {\small $\begin{pmatrix} 2&1\\[-0.5ex] 1&1 \end{pmatrix}$}. The matrix $B-I_2=$ {\small $\begin{pmatrix} 1&1\\[-0.5ex] 1&0 \end{pmatrix}$} is invertible, which already implies $\ZZ^2_r/_{\sim_{B}}=\{(0,0)\}$.\\

\fbox{$T=4$}\quad Here we have exactly one pair $\bz, \bar\bz$ of cycles, see Example~\ref{ExT}. The union  $\bz \cup \bar\bz$ is represented by $B=$ {\small $\begin{pmatrix} 3&2\\[-0.5ex] 1&1 \end{pmatrix}$}. We determine the symmetry group of $B$ according to Remark~\ref{gensymm}. The matrix $B$ has the eigenvalue $\lambda=2+\sqrt 3>1$. Hence $K:=\QQ(\lambda)=\QQ(\sqrt 3)$. The fundamental unit of $\cO_K$ equals $\eps=2+\sqrt 3=\lambda$. Therefore, $\cS(B)$ is generated by $B$ and $-I_2$. The matrix $B$ is reversible since $B^\top\in\bar \bz$. From Example~\ref{ExT}, we see that $B^\top=\hat A \bar B\hat A^{-1}$ holds. Hence 
$$R:=\hat W^{-1}\hat A\hat U=\begin{pmatrix} -1&0\\ 1&1 \end{pmatrix}$$
is a reversing symmetry. Thus $\cR(B)$ is generated by $B,-I_2$ and $R$. Now we can determine $(\ZZ_r)^2/_{\sim_{B}}$. We first consider the equivalence relation generated by (i). We have $B-I_2=$ {\small $\begin{pmatrix} 2&2\\[-0.5ex] 1&0 \end{pmatrix}$}. If $r$ is odd, then there is only one equivalence class, which is represented by $(0,0)$ and we are already done. If $r$ is even, then there are two equivalence classes, namely $\{(l_1,l_2)\mid l_1 \mbox{ is even} \}$ represented by $(0,0)$ and $\{(l_1,l_2)\mid l_1 \mbox{ is odd} \}$, which is represented by $(1,0)$. The matrices $B,-I_2$ and $R$ map each of these sets to itself, thus also $\cR(B)$ preserves these sets. Hence the set $(\ZZ_r)^2/_{\sim_{B}}$ consists of one element represented by $(0,0)$, if $r$ is odd, and
of two elements represented by $(0,0)$ and $(1,0)$, if $r$ is even.\\

\fbox{$T=5$}\quad In this case, $H^0_{\rm red}$ consists of four elements and we have one pair $\bz,\bar\bz$ of cycles, where $\bz$ is represented by $B=$ {\small $\begin{pmatrix} 4&3\\[-0.5ex] 1&1 \end{pmatrix}$}. The eigenvalue $\lambda>1$ equals $\lambda=\frac12(5+\sqrt{21})$. As in the last example, the fundamental unit of $\cO_K$ for  $K:=\QQ(\lambda)=\QQ(\sqrt {21})$ equals $\lambda$, thus $\cS(B)$ is generated by $B$ and $-I_2$. Moreover, also here $R=\hat W^{-1}\hat A\hat U$ is a reversing symmetry and $\cR(B)$ is generated by $B, -I_2$ and $R$. Let us determine $(\ZZ_r)^2/_{\sim_{B}}$. As above, we first consider the equivalence relation generated by (i). Here we have $B-I_2=$ {\small $\begin{pmatrix} 3&3\\[-0.5ex] 1&0 \end{pmatrix}$}.  Therefore, if $3\nmid r$, then there is only one equivalence class, which is represented by $(0,0)$ and we are already done. If $3\mid r$, then we have three equivalence classes, namely,  $O_j:=\{(l_1,l_2)\mid l_1 \equiv j \mbox{ mod } 3\}$ for $j=0,1,2$.  Obviously, $O_j$ is represented by $(j,0)$, $j=0,1,2$. The generators $B,-I_2,R$ of $\cR(B)$ act in the following way on these sets. The matrix $B$ preserves each $O_j$, while $-I_2$ and $R$ preserve $O_0$ and interchange $O_1$ and $O_2$. Consequently, the set $(\ZZ_r)^2/_{\sim_{B}}$ consists of one element represented by $(0,0)$, if $3\nmid r$, and
of two elements represented by $(0,0)$ and $(1,0)$, if $3\mid r$.\\

\fbox{$T=6$}\quad In this case, $H^0_{\rm red}$ consists of six elements and decomposes into three cycles. We have one pair $\bz_1$, $\bar \bz_1=\bz_1^\top$, and a single a cycle $\bz_2=\bar \bz_2=\bz_2^\top$. The pair $\bz_1\cup \bar \bz_1$ is represented by $B_1=$ {\small $\begin{pmatrix} 5&4\\[-0.5ex] 1&1 \end{pmatrix}$}, and $\bz_2$ is represented by $B_2=$ {\small $\begin{pmatrix} 5&2\\[-0.5ex] 2&1 \end{pmatrix}$}. The matrices $B_1$ and $B_2$ have the eigenvalue $\lambda=3+2\sqrt 2$. Furthermore, $K:=\QQ(\lambda)=\QQ(\sqrt {2})$ and $\eps=1+\sqrt 2=\frac12(\lambda-1)$ is the fundamental unit of $\cO_K$. 

Let us first consider $B_1$ and proceed as in Remark~\ref{gensymm}. The matrix $\frac12(B_1-I_2)$ corresponding to $\eps$ is not an integer matrix but its square is integral and equals $B_1$. Hence $\cS(B_1)$ is generated by $B_1$ and $-I_2$. As in the above examples, $R$ is a reversing symmetry and $\cR(B_1)$ is generated by $B_1, -I_2$ and $R$.  In order to apply (i), we need $B_1-I_2=$ {\small $\begin{pmatrix} 4&4\\[-0.5ex] 1&0\end{pmatrix}$}. Suppose that $r$ is odd. Then all elements of $(\ZZ_r)^2$ are equivalent to $(0,0)$ with respect to $\sim_{B_1}$, which can be seen from (i). Now let us consider the case where $4\mid r$. Then there are four equivalence classes $O_0, \dots,O_3$ of the equivalence relation $\sim_{\rm (i)}$ defined by (i). They are represented by $(j,0)$, $j=0,\dots,3$. If we identify the set of equivalence classes with $\ZZ_4$ in the obvious way, then we can describe the action of $\cR(B_1)$ on them as follows. While $B_1$ acts as identity, $R$ and $-I_2$ act as $-\Id$. Hence $(\ZZ_r)^2/_{\sim_{B_1}}$ consists of three elements represented by $(0,0)$, $(1,0)$ and $(2,0)$. Finally, suppose that $r\equiv 2$ mod~4. Then $\sim_{\rm (i)}$ has only two equivalence classes $O_0, O_1$, which are represented by $(0,0)$ and $(1,0)$. The reversing symmetry group of $B_1$ preserves each of the two sets $O_0$ and $O_1$. Consequently, $(\ZZ_r)^2/_{\sim_{B_1}}$ consists of two elements represented by $(0,0)$ and $(1,0)$.

Now we turn to $B_2$. The matrix $B':=\frac12(B_2-I_2)=$ {\small $\begin{pmatrix} 2&1\\[-0.5ex] 1&0\end{pmatrix}$} corresponding to $\eps$ is an integer matrix. Hence $\cS(B_2)$ is generated by $-I_2$ and $B'$. Since $B_2^\top=B_2$, the matrix $\hat W$ is a reversing symmetry of $B_2$. We obtain $\cR(B_2)=\la B', -I_2, \hat W\ra$. We have $B_2-I_2=$ {\small $\begin{pmatrix} 4&2\\[-0.5ex] 2&0\end{pmatrix}$}. If $r$ is odd, then (i) implies that all elements of $(\ZZ_r)^2$ are equivalent to $(0,0)$ with respect to $\sim_{B_2}$. If $r$ is even, then there are four equivalence classes with respect to $\sim_{\rm (i)}$. These are $O_{i,j}$, $i,j\in\{0,1\}$, represented by $(i,j)$. The generator $-I_2$ of $\cR(B_2)$ preserves each of these sets. The remaining generators $B'$ and $\hat W$ preserve $O_{0,0}$ and $O_{1,1}$ and interchange $O_{1,0}$ and $O_{0,1}$. Hence, $(\ZZ_r)^2/_{\sim_{B_2}}$ consists of three elements represented by $(0,0)$, $(1,0)$ and $(1,1)$.\\

\fbox{$T=7$}\quad Also here $H^0_{\rm red}$ consists of six elements and decomposes into three cycles. We have one pair $\bz_1$, $\bar \bz_1=\bz_1^\top$, and a single a cycle $\bz_2=\bar \bz_2=\bz_2^\top$. The pair $\bz_1\cup \bar \bz_1$ is represented by $B_1=$ {\small $\begin{pmatrix} 6&5\\[-0.5ex] 1&1 \end{pmatrix}$}, and $\bz_2$ is represented by $B_2=$ {\small $\begin{pmatrix} 5&3\\[-0.5ex] 3&2 \end{pmatrix}$}. The matrices $B_1$ and $B_2$ have the eigenvalue $\lambda=\frac12(7+3\sqrt 5)$. Furthermore, $\QQ(\lambda)=\QQ(\sqrt {5})$ has fundamental root $\eps=\frac12(1+\sqrt 5)=\frac13(\lambda-2)$. 

We first consider $B_1$. Then the matrices corresponding to $\eps$ and $\eps^2$ are not integral and the one that corresponds to $\eps^3$ equals $B_1$. Thus, $\cS(B_1)$ is generated by $B_2$ and $-I_2$. Because of the special structure of $B_1$, which we already met above, $R$ is a reversing symmetry and $\cR(B_1)$ is generated by $B_1, -I_2$ and $R$. As for (i), we have $B_1-I_2=$ {\small $\begin{pmatrix} 5&5\\[-0.5ex] 1&0\end{pmatrix}$}. In particular, if $5\nmid r$, then all elements of $(\ZZ_r)^2$ are equivalent to $(0,0)$ with respect to $\sim_{B_1}$. If $5\mid r$, then there are five equivalence classes with respect to $\sim_{\rm (i)}$. They are represented by $(j,0)$, $j=0,\dots,4$. Using the map $(j,0)\mapsto j\in\ZZ_5$, we identify the set of equivalence classes with $\ZZ_5$ in order to describe the action of $\cR(B_1)$. The matrix $B_1$ acts as identity, $R$ and $-I_2$ act as $-\Id$. Thus $(\ZZ_r)^2/_{\sim_{B_1}}$ consists of three elements represented by $(0,0)$, $(1,0)$ and $(2,0)$.

Finally, we turn to $B_2$. The matrix $B':=\frac13(B_2-2I_2)=$ {\small $\begin{pmatrix} 1&1\\[-0.5ex] 1&0\end{pmatrix}$} corresponding to $\eps$ is an integer matrix. Hence $\cS(B_2)$ is generated by $-I_2$ and $B'$. We again use that $B_2$ is symmetric and conclude that $\hat W$ is a reversing symmetry of $B_2$, which implies $\cR(B_2)=\la B', -I_2,\hat W\ra$. The image of the map $B_2-I_2=$ {\small $\begin{pmatrix} 4&3\\[-0.5ex] 3&1\end{pmatrix}$} equals $\Span_{\Bbb Z}\{(0,5), (1,2)\}$. Hence all elements of $\ZZ_r^2$ are equivalent to $(0,0)$ with respect to $\sim_{B_2}$ if $5\nmid r$. Now suppose that $5\mid r$. Then there are five equivalence classes with respect to $\sim_{\rm (i)}$, which are respresented by $(0,j)$, $j=0,\dots,4$. As above, we identify the set of equivalence classes with $\ZZ_5$ in the obvious way. Then the generators of $\cR(B_2)$ act by multiplication on this set. Indeed, the matrix $B'$ acts by multiplication by $-2$, $-I_2$  by  $-1$, 
and $\hat W$ by $2$. Consequently, $(\ZZ_r)^2/_{\sim_{B_2}}$ consists of two elements represented by $(0,0)$, $(0,1)$. \\

We summarise the results in the following table. For each $T\le 7$, it contains a list of data $B,r,l$, where the possible values for $l$ depend on divisibility properties of $r$. The groups $\dOsc(B,l)$ for these data constitute a complete system of representatives for the isomorphism classes of discrete oscillator groups with $T\le7$.

\begin{center}
{\renewcommand{\arraystretch}{1.5}
\begin{tabular}{|c|c|l|l|} \hline	
$T$		& $B$ 			& \multicolumn{2}{|c|}{$l$} 		\\
\hline 
\hline
3&{\small $\begin{pmatrix} 2&1\\[-1ex] 1&1 \end{pmatrix}$}& \multicolumn{2}{|c|}{$(0,0)$}\\
\hline
{\multirow{2}{*}{4}}& {\multirow{2}{*}{\small{$\begin{pmatrix} 3&2\\[-1ex] 1&1 \end{pmatrix}$}}}& $r\equiv 0\ (2)$ & $(0,0),\ (1,0)$\\
\cline{3-4} 
&&$r\equiv 1\ (2)$ & $(0,0)$\\
\hline
{\multirow{2}{*}{5}}& {\multirow{2}{*}{\small{$\begin{pmatrix} 4&3\\[-1ex] 1&1 \end{pmatrix}$}}}& $r\equiv 0\ (3)$ & $(0,0),\ (1,0)$\\
\cline{3-4} 
&&$r\not\equiv 0\ (3)$ & $(0,0)$\\
\hline
{\multirow{2}{*}{6}}& {\multirow{2}{*}{{\small$\begin{pmatrix} 5&4\\[-1ex] 1&1 \end{pmatrix}$}}}& $r\equiv 0 \ (4)$ & $(0,0),\ (1,0),\ (2,0)$\\
\cline{3-4} 
&&$r\equiv 2\ (4)$ & $(0,0),\ (1,0)$\\
\cline{3-4} 
&&$r\equiv 1\ (2)$ & $(0,0)$\\
\cline{2-4}
{\multirow{2}{*}{}}& {\multirow{2}{*}{\small{$\begin{pmatrix} 5&2\\[-1ex] 2&1 \end{pmatrix}$}}}& $r\equiv 0 \ (2)$ & $(0,0),\ (1,0),\ (1,1)$\\
\cline{3-4} 
&&$r\equiv 1\ (2)$ & $(0,0)$\\
\hline
{\multirow{2}{*}{7}}& {\multirow{2}{*}{{\small$\begin{pmatrix} 6&5\\[-1ex] 1&1 \end{pmatrix}$}}}& $r\equiv 0 \ (5)$ & $(0,0),\ (1,0),\ (2,0)$\\
\cline{3-4} 
&&$r\not\equiv 0\ (5)$ & $(0,0)$\\
\cline{2-4}
{\multirow{2}{*}{}}& {\multirow{2}{*}{\small{$\begin{pmatrix} 5&3\\[-1ex] 3&2 \end{pmatrix}$}}}& $r\equiv 0 \ (5)$ & $(0,0),\ (0,1)$\\
\cline{3-4} 
&&$r\not\equiv 0\ (5)$ & $(0,0)$\\
\hline
\end{tabular}}
\end{center}
\section{Commensurability}
\begin{de}
Two groups $G_1$ and $G_2$ are said to be abstractly commensurable if there are subgroups $H_1 \subset G_1$ and $H_2 \subset G_2$ of finite index such that $H_1$ is isomorphic to $H_2$. Two subgroups $G_{1}$ and $G_{2}$ of the same group $G$ are said to be commensurable if their intersection $G_{1}\cap G_{2}$ is of finite index in both of them. We say that $G_{1}$ and $G_{2}$ are weakly commensurable if there is an automorphism $F$ of $G$ such that $G_{1}$ and $F(G_{2})$ are commensurable.
\end{de}
Obviously, two subgroups $G_{1}$ and $G_{2}$ of a group $G$ are weakly commensurable if and only if there are finite index subgroups $H_1\subset G_{1}$ and $H_2\subset G_{2}$ and an automorphism $F$ of $G$ such that $F(H_2)=H_1$. Finite index subgroups of lattices are again lattices. Hence Theorem~\ref{isomequiv} implies that two lattices in $\Osc$ are weakly commensurable if and only if the underlying discrete groups are abstractly commensurable.

\begin{theo}\label{Tcomm}
For an element $B\in \SL(2,\ZZ)$ with $\tr B>2$, let $\lambda(B)$ be any of the two eigenvalues of $B$. The map 
$${\cG}\longrightarrow \cK,\quad \dOsc(B,l)\longmapsto \QQ(\lambda(B)) $$
induces a bijection from the set of abstract commensurability classes of discrete split oscillator groups to the set $\cK$ of real quadratic fields.
\end{theo}
Note that $\QQ(\lambda(B))=\QQ(\sqrt{T^2-4})$ for every $B\in \bB_T$, $T>2$.

In order to prove Theorem~\ref{Tcomm} we need the following lemmata.
\begin{lm} {\rm \cite{Ha}} \label{finite}
Let $G$ be a finitely generated group. For a fixed $n \in \NN$, $G$ has a finite number of subgroups of index $n$.
\end{lm}

We define an equivalence relation on $\GL(2,\ZZ)$ by $B_1\sim B_2$ if and only if there exist integers $m,n \in \ZZ_{\not=0}$ and a matrix $M \in \operatorname{GL(2, \QQ)}$ such that $B_1^m = MB_2^{n}M^{-1}$.
\begin{lm}\label{equiv}
Let $B_1,B_2 \in \bB$ have eigenvalues $\lambda_1^{\pm1}$ and  $\lambda_2^{\pm1}$, respectively. 
Then $B_1 \sim B_2$ if and only if $\QQ(\lambda_1) = \QQ(\lambda_2)$.
\end{lm}
\proof The characteristic polynomial $p_B$ for any $B\in\bB$ is irreducible over $\QQ$. Hence $B$ is conjugate over $\QQ$ to the companion matrix of $p_B(x)=x^2-Tx+1$, where $T=\tr(B)$ (Frobenius normal form). Thus two matrices $B, B'\in\bB$ are conjugate in $\GL(2,\QQ)$ if and only if they have the same trace. Consequently, $B_1\sim B_2$ holds for $B_1,B_2\in\bB$ if and only if  there exist $n,m \in \ZZ_{\not=0}$ such that $\lambda_1^m + 1/\lambda_1^m = \tr(B_1^m) = \tr(B_2^n) = \lambda_2^n + 1/\lambda_2^n$, which is equivalent to  $\lambda_1^m = \lambda_2^n$ or $\lambda_1^m = \lambda_2^{-n}.$

Suppose that $B_1\sim B_2$. Then the above considerations imply that $\lambda_1^m$ is in $\QQ(\lambda_1)\cap\QQ(\lambda_2)$. This shows that the latter intersection is not equal to $\QQ$. Indeed, $\lambda_1^m$ is a root of the characteristic polynomial of $B_1^m$, which has integer coefficients and is irreducible over $\ZZ$ thus also over $\QQ$. Hence $\lambda_1^m$ is irrational. Consequently, $\QQ(\lambda_1)=\QQ(\lambda_2)$.

Conversely, suppose that $\QQ(\lambda_1)=\QQ(\lambda_2)=:K$. Since $\lambda_1$ and $\lambda_2$ are units in $\cO_K$, they are positive or negative powers of the fundamental unit. Hence there are integers $m,n\in\ZZ_{\not=0}$ such that $\lambda_1^m=\lambda_2^n$. The above argument now implies $B_1\sim B_2$.
\qed
\begin{lm}\label{concomm}
Two discrete split oscillator groups ${\rm Osc}^{r_1}(B_1,l_1)$ and ${\rm Osc}^{r_2}(B_2,l_2)$ are abstractly commensurable if and only if $B_1 \sim B_2$.
\end{lm}
\proof In this proof `commensurable' always means `abstractly commensurable'. We introduce the notation $\Gamma(B):={\rm Osc}^1(B,0)$. First we claim that every discrete split oscillator group $\dOsc(B,l)$ contains a finite index subgroup that is isomorphic to $\Gamma(B^k)$ for some $k\in\NN$. To show this, we consider the discrete Heisenberg group $H^r(\ZZ)$, $r\in\NN$, $r>0$, as defined in~(\ref{dHeis}). Its subgroup generated by $\alpha, \beta, \gamma^r$ is isomorphic to $H^1(\ZZ)$. Of course, the action of $\delta\in\dOsc(B,l)$ does not necessarily preserve this subgroup. However, by Lemma~\ref{finite} some power of $\delta$ does. Let $\delta':=\delta^k$ be this power.  Then the group generated by $\alpha, \beta, \gamma^r$ and $\delta'$ is a finite index subgroup of $\dOsc(B,l)$. It is isomorphic to $H^1(\ZZ)\rtimes \ZZ\delta'$, hence to $\Gamma(B^k)$, which proves the claim.

Next assume that $\Gamma(B)$ and $\Gamma(B')$ are commensurable. Let $\Gamma$ be a finite index subgroup of $\Gamma(B)$ that can be embedded as a finite index subgroup into  $\Gamma(B')$. Then $\Gamma$ is also a discrete split oscillator group. Indeed, $\Gamma(B)$ can be embedded as a lattice into $\Osc$. Hence, the image of $\Gamma$ is also a lattice in $\Osc$, thus isomorphic to a discrete split oscillator group. More exactly, $\Gamma$ is isomorphic to $\dOsc(B^k,l)$ for suitable $r,k\in\NN_{>0}$ and $l\in(\ZZ_r)^2$. On the other hand, $\Gamma$ is isomorphic to a finite index subgroup of $\Gamma(B')$. Thus the same reasoning implies that $\Gamma$ is isomorphic to a group ${\rm Osc}^{r'}(B'^j,l')$ for suitable $r',j\in\NN_{>0}$ and $l'\in\ZZ_r^2$. Proposition~\ref{pre} now implies that $B^k$ is conjugate to $B'^j$ or to $(B')^{-j}$.

Now we can prove the lemma. Suppose that ${\rm Osc}^{r_1}(B_1,l_1)$ and ${\rm Osc}^{r_2}(B_2,l_2)$ are commensurable. Then $\Gamma(B_1^p)$ and $\Gamma(B_2^q)$ are commensurable for suitable $p,q\in\NN_{>0}$ by the above consideration. Hence there exist numbers $k,j\in\NN_{>0}$ such that $B_1^{pk}$  is conjugate to $B_2^{jq}$ or to $B_2^{-jq}$. This implies $B_1\sim B_2$.

Conversely, suppose that $B_1\sim B_2$ holds. Choose $m,n\in\ZZ_{\not=0}$ such that $B_1^m$ and $B_2^n$ are conjugate. Furthermore, choose $p,q\in\NN_{>0}$ such that ${\rm Osc}^{r_1}(B_1,l_1)$ is commensurable with $\Gamma(B_1^p)$ and ${\rm Osc}^{r_2}(B_2,l_2)$ is commensurable with $\Gamma(B_2^q)$. We have to show that $\Gamma(B_1^p)$ and $\Gamma(B_2^q)$ are commensurable. We use that $\Gamma(B_1^p)$ and $\Gamma_1:=\Gamma((B_1^p)^{mq})$ are commensurable and that $\Gamma(B_2^q)$ and $\Gamma_2:=\Gamma((B_2^q)^{np})$ are commensurable. Since $B_1^{mpq}$ and $B_2^{npq}$ are conjugate, $\Gamma_1$ and $\Gamma_2$ are isomorphic by Proposition~\ref{pre}.
\qed

{\sl Proof of Theorem~\ref{Tcomm}.}\ Lemmata~\ref{equiv} and~\ref{concomm} show that the map $\cG\rightarrow \cK$ indeed descends to a map from the set of abstract commensurability classes of discrete split oscillator groups to the set $\cK$ and that this map is injective. It remains to show that it is also surjective. Let $K$ be a real quadratic field. We choose a unit $\lambda>1$ of norm 1 in $\cO_K$. Then $K=\QQ(\lambda)$. The minimal polynomial of $\lambda$ is of the form $p(x)=x^2+a x+1$, $a\in\ZZ$. The companion matrix $C_p$ of $p$ is in $\SL(2,\ZZ)$ and $\lambda$ is an eigenvalue of $C_p$. Consequently, the isomorphism class of $\Gamma(C_p)$ is a preimage of $K$.
\qed
\section{Compact Clifford-Klein forms of the symmetric space $\Osc$}

The aim of this subsection is to relate our results to those obtained by Maeta in \cite{Ma}.
Maeta considers compact quotients of solvable pseudo-Riemannian symmetric spaces by discrete subgroups of the transvection group. He proved that for signature $(2,2)$ such quotients exist only for two solvable symmetric spaces. One of these spaces is the split oscillator group, which becomes a symmetric space if we endow it with the biinvariant metric defined in Subsection~\ref{S21}.

Of course, each lattice $L$ of $\Osc$ gives rise to a compact quotient. We will see that these quotients cannot be obtained as a Clifford-Klein form if we consider the symmetric space $\Osc$ as a quotient $\hat G/\hat G_+$ of its transvection group $\hat G$ by the stabiliser $\hat G_+$. Recall that the isometry group of a pseudo-Riemannian symmetric space can have a larger dimension than the transvection group. Indeed, for the symmetric space $\Osc$, the transvection group $\hat G$ is properly contained in the identity component of the isometry group, which we want do denote by $G$. In particular, we will see that the subgroup of $G$ that consists of left-translations by group elements of $\Osc$ is not contained in $\hat G$. Consequently, we have to consider $\Osc$ as a homogeneous space $G/G_+$ in order to understand $L\backslash \Osc$ as a Clifford-Klein form. Here $G_+\subset G$ denotes the stabiliser of the base point. In particular,  the compact quotients of the symmetric space $\Osc$ considered in~\cite{Ma} are not of the form $L\backslash \Osc$.

Let us explain this in more detail. We use the two-form $\omega$ and the map $A\in\SO^+(1,1)$ defined in Section~\ref{S21} in order to define a 6-dimensional  generalised oscillator group. As usual, we identify $\RR^2\cong \CC$ and consider the 5-dimensional Heisenberg group $H_2=\RR\times\CC\times\CC$ with group multiplication
$$(z,\xi_1,\xi_2)\cdot(\hat z,\hat\xi_1,\hat\xi_2)= \big(z+\hat z +\textstyle \frac12 (\omega(\xi_1,\hat\xi_1)-\omega(\xi_2,\hat\xi_2)),\xi_1+\hat \xi_1, \xi_2+\hat\xi_2\big).$$ 
We define an action of $\RR$ on $H_2$ by
$$t.(z,\xi_1,\xi_2)=(z, e^{tA}\xi_1, e^{tA} \xi_2).$$
Using this action we obtain a semidirect product
$${\rm Osc}_{2,0}(1,1):=H_2\rtimes \RR.$$
This notation for the generalised oscillator group was introduced  in \cite{KO}. 

Since we want to distinguish between the symmetric space $\Osc$ and the group $\Osc$, we will denote the symmetric space by $X$ and the group by $Q$. 

\begin{pr} The transvection group $\hat G$ of $X$ is isomorphic to ${\rm Osc}_{2,0}(1,1)$. The unity component $G$ of the isometry group of $X$ is isomorphic to ${\rm Osc}_{2,0}(1,1)\rtimes \RR$, where $\RR$ acts on ${\rm Osc}_{2,0}(1,1)$ by
\begin{equation}\label{sacts}
s.(z,\xi_1,\xi_2,t)=(z, e^{s\phi}(\xi_1,\xi_2),t),
\end{equation}
and $\phi(\xi_1,\xi_2)=(\xi_2,\xi_1)$. 
The stabilisers of the base point are equal to $\hat G_+=\{ (0, iy_1, iy_2,0)\mid y_1,y_2\in\RR\}$ and $G_+=\hat G_+\rtimes \RR$, respectively.

The map 
\begin{equation}\label{Q}
\iota:\ Q\hookrightarrow G, \quad (2z,2\xi,2t)\longmapsto (z,\xi,  A\xi,t,  t)
\end{equation}
is a homomorphism. Its image equals the group of left translations on $X$.
\end{pr}
\proof
The action of $Q\times Q$ on $X$ defined by $(q_1,q_2)\cdot x=q_1xq_2^{-1}$  is isometric since the metric on $X$ is biinvariant. The kernel of this action equals $\{(z,z)\in Q\times Q\mid z\in Z(Q)\}\cong Z(Q)$. Hence $I:=(Q\times Q)/Z(Q)$ is a subgroup of the isometry group of $X$. This subgroup acts transitively on $X$. It is invariant under the conjugation by the geodesic reflection of $X$ at the identity, which we denote by $\theta$. In particular, it contains the transvection group $\hat G$. The Lie algebra $\fri$ of $I$ equals $(\fq\oplus\fq)/\fz(\fq)$, where $\fq$ denotes the Lie algebra of $\fq$. Now we consider the eigenspace decomposition of $\fri$ with respect to the differential of $\theta$. The (-1)-eigenspace $\fri_-$ equals the anti-diagonal $\{(X,-X)\mid X\in\fq\}$. The (+1)-eigenspace $\fri_+$ equals $\{(X,X)\mid X\in \fq\}/\fz(\fq).$ Next we turn to the Lie algebra $\hat \fg$ of the transvection group $\hat G\subset I$. While its $(-1)$-eigenspace $\hat\fg_-$ equals $\fri_-$, the (+1)-eigenspace equals $\hat \fg_+=\{(X,X)\mid X\in[\fq,\fq]\}/(\fz(\fq)\cap [\fq,\fq])\subset \fri_+.$ A direct calculation shows that $\hat \fg=\hat\fg_+\oplus\hat\fg_-$ can be identified with ${\rm Osc}_{2,0}(1,1)$ such that 
$$\hat\fg_+=\{(0,iy_1,iy_2,0)\mid y_1,y_2\in \RR\},\quad \hat\fg_-=\{(z,x_1,x_2,t)\mid x_1,x_2,z,t\in\RR\}$$
under this identification. 

Next we want to show that the unity component $G$ of the isometry group is equal to $I$. Note that the elements of $\fg_+$ have to act as derivations on $\hat \fg$. Moreover, they have to preserve the subspaces $\hat \fg_+$ and $\hat \fg_-$ and the scalar product on $\hat\fg_-\cong T_eX$. It is easy to see that this implies $\dim\fg_+\le 3$. Since $\fri_+\subset\fg_+$ and $\dim \fri_+=3$, we get $\fg_+=\fri_+$, which proves the assertion. Since $\dim \hat\fg_+=2$, we obtain $G=\hat G\rtimes \RR$. Under the identification of $\hat G$ with ${\rm Osc}_{2,0}(1,1)$ the action of $\RR$ on $\hat G$ is given by (\ref{sacts}).

Finally, the subgroup of left translations equals $Q\times\{e\}\subset I$. With the identification of $I$ and ${\rm Osc}_{2,0}(1,1)\rtimes\RR$, this subgroup equals $\{(z,\xi,A\xi,t,t)\mid \xi\in\CC, z,t\in\RR\}$. It can be checked directly that the map in the proposition that maps $Q$ to this subgroup is a homomorphism.
\qed

Finally we want to show that quotients by lattices are straight in the sense of~\cite{KO}. Moreover, we will see that there are also non-straight quotients. This is analogous to the case of the (ordinary) oscillator group, which has been considered in~\cite{KO}. Let $L$ be a discrete subgroup of $G$ acting freely and properly on $X$ such that $\Gamma\backslash X$ is compact. We consider the projection of $\Gamma$ to the $\RR$-factor of $\hat G=H_2\rtimes\RR\subset G$, i.e., the subgroup $p(\Gamma)\subset \RR$ for  
$$p:G=H_2\rtimes (\RR\oplus \RR)\longrightarrow \RR,\quad (z,\xi_1,\xi_2,t,s)\longmapsto t.$$
The quotient $\Gamma\backslash X$ is called straight if $p(\Gamma)\subset \RR$ is discrete. Otherwise it is called non-straight. 
\begin{lm}
The set $N:=\{(z,\xi,-A\xi,t,t)\mid \xi\in\CC,\, z,t\in\RR\}\subset G$ is a subgroup of $G$. It is isomorphic to the direct product $H\times \RR$ of the three-dimensional Heisenberg group $H$ and $\RR$. The map
\begin{equation}\label{N}
N\longrightarrow X= G/G_+,\quad n\longmapsto nG_+
\end{equation}
is a diffeomorphism. The action of $N$ on $X$  corresponds to the left translation on $N$ under this diffeomorphism. 
\end{lm}
\proof 
The map $\phi+\diag(A,A)$ acts trivially on $\{(\xi,-A\xi)\mid \xi\in\CC\}$ since
$$(\phi+\diag(A,A))(\xi,-A\xi)=(-A\xi,\xi)+(A\xi,-A^2\xi)=0.$$
Hence $\fn:=\{ (z,\xi,-A\xi,t,t)\mid  \xi\in\CC,\, z,t\in\RR\}\subset\fg$ is a subalgebra. The set $N\subset G$ is the analytic subgroup corresponding to $\fn$, which proves the first assertion. The second one follows from 
$\omega(\xi_1,\xi_2)-\omega(-A\xi_1,-A\xi_2)=2\omega(\xi_1,\xi_2).$

Since $\fg_+=\{(0,iy_1,iy_2,0,s)\mid y_1,y_2,s\in\RR\}$, we obtain that $\fg=\fn\oplus\fg_+$ is a direct sum of subalgebras. Hence, $G=N\cdot G_+$. Moreover, $G_+\cap N=\{e\}$, since $(z,\xi,-A\xi,t,t)=(0,iy_1,iy_2,0,s)$ implies $z=\xi=y_1=y_2=t=s=0$. This yields the remaining assertions. \qedohne
%
%
%
\begin{pr} Quotients of $X=\Osc$ by lattices are straight. The symmetric space $X$ also admits non-straight quotients by discrete subgroups of the isometry group.
\end{pr}
\proof 
Let $L$ be a lattice in $Q=\Osc$. In Section~\ref{S61}, we have seen that the projection of $L$ to the $\RR$-factor of $\Osc=H\rtimes\RR$ is discrete. Hence also $p(\iota(L))$ is discrete, where $\iota: Q\hookrightarrow G$ is the embedding defined in (\ref{Q}). This proves the first assertion.

In order to construct non-straight lattices, we use the diffeomeorphism defined by (\ref{N}). We see that every lattice in $N$ defines a compact quotient of $X$. Recall that $N\cong H\times\RR$, where the projection to the $\RR$-factor corresponds to $p|_N$ under this isomorphism. Thus it suffices to construct lattices  $\Gamma$ in $H\times \RR$ whose projection to $\RR$ is dense in $\RR$.  Let $\Gamma_0$  be a lattice of the Heisenberg group $H$, and let $\ph:H \to \RR$ be a group homomorphism such that $\ph(\Gamma_0)$ is not contained in $\QQ$. Then
$\Gamma:=\{(\gamma_0,\ph(\gamma_0)+k)\mid \gamma_0\in\Gamma_0,\, k\in\ZZ\}$ is such a lattice.
\qed
\begin{re} {\rm Maeta \cite{Ma} constructed compact quotients of $X=\Osc$ by discrete subgroups of the transvection group $\hat G$. He uses that  $\hat G=L_{\Id,0}\cdot \hat G_+$, where
$$L_{\Id,0}= \{(z,x+ix,y-iy,t)\mid x,y,z,t\in\RR\}.$$
Obviously, $L_{\Id,0}\cap\hat G_+=\{e\}$. Note that $L_{\Id,0}\cong \RR^3\rtimes\RR$, where $t\in\RR$ acts on $\RR^3$ by $e^{tM}$ for $M=\diag(0,1,-1)$. Analogously to the above considerations, it suffices to construct lattices in $L_{\Id,0}$. Such lattices can be obtained in the following way. Let $r$ be a positive real number such that $r+r^{-1}$ is an integer and put $t_0:=\ln r$. Then the characteristic polynomial of $e^{t_0M}$ has integer coefficients. Hence there exists a lattice $\Gamma_0$ in $\RR^3$ that is preserved by $e^{t_0M}$. Consequently, $\Gamma_0\rtimes t_0\ZZ$ is a lattice in $L_{\Id,0}$. 

Compact quotients of $X=\Osc$ by lattices of $L_{\Id, 0}$ are straight. Indeed, if $\Gamma$ is a lattice in $L_{\Id,0}$, then the projection of $\Gamma\subset\RR^3\rtimes\RR$ to the $\RR$-factor is discrete since the action of $\RR$ on $\RR^3$ is non-trivial, compare the argument for $\Osc$ in Section~\ref{S61}.  This projection is equal to the projection $p$ on $L_{\Id,0}$. Hence, $\Gamma\backslash X$ is straight.
}
\end{re}
\appendix
\section{Appendix: On conjugacy classes of integer matrices}
Recall that we defined $\bB_T=\{B\in\SL(2,\ZZ)\mid \tr B=T\}$ for $T>2$. The group $\SL(2,\ZZ)$ acts by conjugation on this set. Obviously, this is in fact an action of ${\rm PSL}(2,\ZZ)$. We want to describe the orbits of this action. There are various known methods to do this. Here we want to present the one developed by Aicardi~\cite{A}, which we adapt to our situation.

Recall that $\SL(2,\ZZ)$ is generated by the elements
$$\hat A:=\begin{pmatrix}1&1\\0&1 \end{pmatrix},\quad \hat B:=\begin{pmatrix}1&0\\1&1 \end{pmatrix},\quad  \hat R:=\begin{pmatrix}0&1\\-1&0 \end{pmatrix}.$$
\begin{de}
A cycle of length $t>1$ in $\bB_T$ is a cyclic sequence $[B_1, B_2,\dots, B_t]$ of distinct matrices such that $B_i = M_{i-1}B_{i-1}M_{i-1}^{-1}$ $(i = 2,\dots,t)$ and $B_1 = M_t B_t M_t ^{-1}$, where each of the matrices
$M_1,\dots,M_t$ equals $\hat A$ or $\hat B$. The matrices $M_1,\dots,M_t$ will be called cycle operators.
\end{de}
 In the following, we consider the index $i$ of the elements in a cycle of length $t$ as an element of $\ZZ_t=\ZZ/t\ZZ$, i.e., we put $i+1=1$ if $i=t$ and $i-1=t$ if $i=1$.
 
Now let us consider the subset
$$H^0=\left\{ \begin{pmatrix}a&b\\c&d\end{pmatrix} \in\bB_T \mid a,b,c,d>0\right\}$$
of $\bB_T$. A cycle $[B_1, B_2,\dots, B_t]$ in $\bB_T$ is called cycle in $H^0$, if $B_j\in H^0$ for all $j=1,\dots, t$.
\begin{pr}\label{c}
 Each ${\rm PSL}(2,\ZZ)$-orbit in $\bB_T$ contains exactly one cycle in $H^0$.
\end{pr}
\proof We show that the assertion directly follows from \cite{A}, Theorem  4.10.
We define a homomorphism $\ph: \SL(2,\ZZ) \longrightarrow \SO(2,1)$
by $\ph(\hat R)=\diag(-1,-1,1)$ and
$$\ph(\hat A)=\begin{pmatrix}1&1&1\\-1&1/2&-1/2\\1&1/2&3/2 \end{pmatrix},\quad \ph(\hat B)=\begin{pmatrix}1&-1&1\\1&1/2&1/2\\1&-1/2&3/2 \end{pmatrix}.$$
The kernel of $\ph$ equals $\{I_2,-I_2\}$, thus the image of $\ph$ is isomorphic to ${\rm PSL}(2,\ZZ)$. In \cite{A}, this image is denoted by ${\cal T}$.

We consider the subset $\bH_T$ of the hyperboloid $x^2+y^2-z^2=\Delta_T:=T^2-4$ defined by 
$$\bH_T:=\{[K,D,S]\in \ZZ^3\mid K^2+D^2-S^2=\Delta_T,\ D\equiv S \mbox{ mod } 2 \}.$$ 
In the language of \cite{A}, this is the set of good points in $H_\Delta$ for $\Delta=\Delta_T$.  The map
$$\phi:\ \bB_T \longrightarrow \bH_T,\quad  \begin{pmatrix}a&b\\c&d \end{pmatrix}\longmapsto [K,D,S]:=[a-d,b+c,-b+c]$$
is a bijection and
$$\phi(M B M^{-1})=\ph(M)(\phi(B))$$
holds for all $M\in\SL(2,\ZZ)$. Hence, we can identify the action of ${\rm PSL}(2,\ZZ)$ on $\bB_T$ with the action of ${\cal T}$ on $\bH_T$. The image of $H^0$ under $\phi$ is the set
$$\{ [K,D,S]\in\bH_T\mid |S|<D,\ D>0\},$$
which is also denoted by $H^0$ in \cite{A}. Note that $\Delta_T$ is never a square for $T>2$. Thus \cite{A}, Theorem 4.10, applies and gives the assertion. 
\qed

Proposition \ref{c} shows that the finite set $H^0$ decomposes into cycles and every ${\rm PSL}(2,\ZZ)$-orbit in $\bB_T$ contains exactly one of these. If we choose an element in each of these cycles, we obtain a system of representatives of the conjugacy classes contained in $\bB_T$. The following lemma tells us how the decomposition of $H^0$ can be determined. It directly follows from \cite{A}, Lemma 4.6 and \cite{A}, Figure 9.
\begin{lm}\label{LA}
Let $B$ be in $H^0$. Then either $\hat A B\hat A^{-1}$ or $\hat B B\hat B^{-1}$ is in $H^0$. Furthermore, either $\hat A^{-1} B\hat A$ or $\hat B^{-1} B\hat B$ is in $H^0$. 
\end{lm}
Hence we can proceed in the following way. We start by choosing a matrix $B\in H^0$. The next element in the cycle of $B$ is $\hat A B\hat A^{-1} $ if this matrix is in $H^0$, otherwise it equals $\hat BB\hat B^{-1} $.  We continue in this way until we reach again $B$. This gives us the first cycle. Then we proceed in the same way with a matrix that is not in this cycle, etc. Since $H^0$ is finite, this procedure finishes and gives the wished decomposition of $H^0$ into cycles. 

In the following we want to modify this method of describing the  ${\rm PSL}(2,\ZZ)$-orbits in $\bB_T$ in order to make it more efficient. Instead of $H^0$  we consider the set
$$H^0_{\rm red}=\left\{ B=\begin{pmatrix}a&b\\c&d\end{pmatrix} \in\bB_T \mid a,b,c,d>0,\ \max\{a,b,c,d\}\in\{a,d\} \right\}$$
\begin{re}\label{Rmax}{\rm
Note that for $B\in H^0_{\rm red}$ the condition $a=\max\{a,b,c,d\}$ is equivalent to $d=\min\{a,b,c,d\}$. Similarly, $d=\max\{a,b,c,d\}$ is equivalent to $a=\min\{a,b,c,d\}$. Indeed, suppose that $a=\max\{a,b,c,d\}$ holds and assume that $d>b$. Then $1=ad-bc>ad-dc=d(a-c)\ge 0$. Hence $a=c=1$. Now $a=\max\{a,b,c,d\}$ implies that also $b=d=1$, which is a contradiction. Thus $d\le b$. Analogously, $d\le c$ holds. Hence, $d=\min\{a,b,c,d\}$. The remaining statements follow in a similar way.
}
\end{re}
\begin{de}
A reduced cycle in $H^0$ is a cyclic sequence $[B_1, B_2,\dots, B_t]$ of distinct matrices $B_1,\dots, B_t\in H^0_{\rm red}$ such that $B_{i+1} = M_{i}B_{i}M_{i}^{-1}$ for all indices $i$, where each matrix $M_i$ is a power of  $\hat A$ or $\hat B$ and where $M_i$ is a power of $\hat A$ if and only if $M_{i+1}$ is a power of $\hat B$. 
Here the index $i$ is again taken modulo $t$ as explained above. 
\end{de}
By definition, a reduced cycle in $H^0$ is the same as a cycle in $H^0_{\rm red}$ in the sense of Subsection~\ref{S81}. 

\begin{figure}
\centering
{\small
\begin{tikzpicture}
\blue \node (l) at (2.9,2) {$\begin{pmatrix} 4&3\\ 1&1 \end{pmatrix}$\ };
\node (r) at (5.2,2) {\ $\begin{pmatrix} 1&3\\ 1&4 \end{pmatrix}$\ };
\draw[->,bend left=45 ] (l) edge node [above] {$\hat B$} (r); 
\draw[->,bend left=45 ] (r) edge node [above] {$\hat A^3$} (l); 
\black
\node (lu) at (3,0) {$\begin{pmatrix} 3&5\\ 1&2 \end{pmatrix}$\ };
\node (ru) at (5.2,0) {\ $\begin{pmatrix} 2&5\\ 1&3 \end{pmatrix}$\ };
\draw[<-,bend left=-45] (lu) edge node [below] {$\hat A$} (ru);
\draw[<-,bend left=-45] (l) edge node [left] {$\hat A$} (lu);
\draw[->,bend left=45] (r) edge node [right] {$\hat A$ } (ru);
\end{tikzpicture}
}
\caption{Reduced cycle.}
\label{fig}
\end{figure}
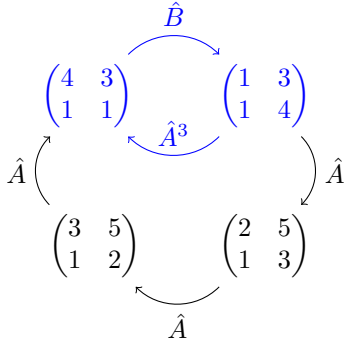

With each cycle $\bz=[B_1,\dots,B_t]$ in $H^0$, we will associate a reduced cycle $\bz_{\rm red}$ by combining all consecutive cycle operators that are equal as shown in Figure~\ref{fig}, where $\bz_{\rm red}$ is displayed in blue. More exactly, we put $\bz_{\rm red}=[B'_1,\dots,B'_r]$, where $\{B'_1,\dots,B'_r\}\subset \{B_1,\dots,B_t\}$ and $B_j$ belongs to $\{B'_1,\dots,B'_r\}$ if and only if  the cycle operators $M_{j-1}$ and $M_j$ are different. The order of the matrices $B'_1,\dots,B'_r$ in $\bz_{\rm red}$ is the same as in $\bz$. 

\begin{pr}\label{P1red}
For each cycle $\bz=[B_0,\dots, B_t]$ in $H^0$, the sequence $\bz_{\rm red}$ is a reduced cycle. It consists of exactly those $B_j\in\bz$ that belong to $H^0_{\rm red}$. For the cycle operators $M'_1,\dots, M'_r$ of $\bz_{\rm red}$, the following holds. The operator $M'_j$ is a power of $\hat A$  if and only if $M'_{j+1}$ is a power of $\hat B$. 
\end{pr}

Before we prove this assertion, we want to show a technical lemma. 
As usual, we use the notation
$$
B=\begin{pmatrix} a&b\\c&d \end{pmatrix},\quad B_i=\begin{pmatrix} a_i&b_i\\c_i&d_i\end{pmatrix}.
$$
Then
\begin{eqnarray}\label{nn1}
\hat A^nB \hat A^{-n} &=&\begin{pmatrix}a+nc& b+n(d-a)-n^2c\\ c &d-nc\end{pmatrix},\\[1ex]
\hat B^nB\hat B^{-n} &=&\begin{pmatrix} a-nb& b\\  c+n(a-d)-n^2b& d+nb\end{pmatrix}. \label{nn2}
\end{eqnarray}
\begin{lm}\label{Lemma}
An element $B_i$ of a cycle $[B_1, \dots, B_t]$ in $H^0$ satisfies $B_i=\hat A B_{i-1} \hat A^{-1}$ and $B_{i+1}=\hat B B_i\hat B^{-1}$ if and only if $a_i=\max\{a_i,b_i,c_i,d_i\}$. Similarly,  $B_i=\hat B B_{i-1} \hat B^{-1}$ and $B_{i+1}=\hat AB_i\hat A^{-1}$ holds if and only if $d_i=\max\{a_i,b_i,c_i,d_i\}$.
\end{lm}
\proof Suppose that $B_i=\hat A B_{i-1} \hat A^{-1}$ and $B_{i+1}=\hat B B_i\hat B^{-1}$. Then  $B_i$, $\hat A^{-1} B_{i} \hat A$ and $\hat BB_i\hat B^{-1}$ are in $H^0$. Therefore all their entries are positive. We use formulas (\ref{nn1}) and (\ref{nn2}) for $n=-1$ and $n=1$, respectively. This yields $a_i>c_i$ and $a_i>b_i$.  Now $\det(B_i)=1$ implies that also $a_i> d_i$ holds. Conversely, suppose that $a_i=\max\{a_i,b_i,c_i,d_i\}$. Then we have $d_i=\min\{a_i,b_i,c_i,d_i\}$ by Remark~\ref{Rmax}. Now (\ref{nn1}) implies that $\hat A B_{i} \hat A^{-1}$ is not in $H^0$. Hence $B_{i+1}=\hat B B_i\hat B^{-1}$. Furthermore, (\ref{nn2}) shows that $\hat B^{-1} B_i\hat B$ is not in $H^0$, which implies $B_{i-1}=\hat A^{-1} B_i\hat A$. The proof of the second assertion is analogous.\qed

{\sl Proof of Proposition~\ref{P1red}.} Let us first note that $\bz_{\rm red}$ is not empty. Indeed, if it were empty, then all cycle operators of $\bz$ would be equal. Hence $B_1=\hat A^nB_1 \hat A^{-n}$ or $B_1=\hat B^nB_1 \hat B^{-n}$ for some $n>0$, which is impossible since $B_1$ is in $H_0$.

By construction, all cycle operators are powers of $\hat A$ or $\hat B$. In order to prove that $\bz_{\rm red}$ is a reduced cycle, it remains to show that all matrices $B'_j$ belong to $H^0_{\rm red}$. However, this is a consequence of Lemma~\ref{Lemma}.  The same Lemma also implies that $\bz_{\rm red}$ contains all $B_j\in\bz$ that belong to $H^0_{\rm red}$. The last assertion of the proposition follows from the construction of $\bz_{\rm red}$.\qed

The following proposition shows that we have a one-to-one correspondence between cycles and reduced cycles in $H^0$. Roughly speaking, the reduced cycles in $H^0$ are exactly the intersections of the cycles in $H^0$ with $H^0_{\rm red}$.
\begin{pr}\label{redunred} The map $R:\bz\mapsto \bz_{\rm red}$ is a bijection from the set of cycles in $H^0$ to the set of reduced cycles in $H^0$.
\end{pr}
\proof  We show the assertion by constructing an inverse $U$ of $R$. Let $\bz'$ be a reduced cycle in $H^0$. All elements of $\bz'$ are conjugate to each other. Hence $\bz'$ is contained in a unique cycle $\bz$ in $H^0$ and we put $U(\bz')=\bz$. Obviously, this is a left inverse of $R$. Let us show that it is also a right inverse. 
Let $\bz'$ be a reduced cycle in $H^0$. We have to show that $\bz_{\rm red}=\bz'$ for $\bz=U(\bz')$.  By construction, $\bz'$ is a reduced cycle contained in $\bz$. Since $\bz_{\rm red}$ contains all elements of $\bz$ that are in $H^0_{\rm red}$, the elements of $\bz'$ are contained in $\bz_{\rm red}$. Assume that $\bz'$ were not equal to $\bz_{\rm red}$. Then there would exist an element $B$ of $\bz'$ whose successor $B'$ in $\bz'$ is different from that in $\bz_{\rm red}$, which we will denote by $B_{\rm red}$. Suppose that $d=\max\{a,b,c,d\}$ holds for $B$. Then $a=\min\{a,b,c,d\}\le b$, thus $\hat B^n B \hat B^{-n}\not\in H^0$ for all $n\ge1$, see (\ref{nn2}). This would imply $B'=\hat A^k B\hat A^{-k}$ and $B_{\rm red}=\hat A^l B\hat A^{-l}$ for some positive numbers $k\not=l$. Let $a_{\rm red},b_{\rm red},c_{\rm red},d_{\rm red}$ denote the entries of $B_{\rm red}$. Since $B_{\rm red}$ is the successor of $B$ in $\bz_{\rm red}$, the matrices $\hat A^n B\hat A^{-n}$, $n=1,\dots,l-1$, would belong to $\bz$ but not to $H^0_{\rm red}$ by Lemma~\ref{Lemma}. Thus $p:=k-l>0$. Then $\hat A^{p} B_{\rm red} \hat A^{-p}=B'$ belongs to $H^0$. Thus $d_{\rm red}> pc_{\rm red}$ by (\ref{nn1}). In particular, $d_{\rm red}$ is not the minimal entry of $B_{\rm red}$. Hence $\min\{a_{\rm red},\dots,d_{\rm red}\}=a_{\rm red}$ since $B_{\rm red}\in H^0_{\rm red}$. On the other hand, the successor of $B_{\rm red}$ in $\bz_{\rm red}$ equals $\hat B^{q} B_{\rm red} \hat B^{-q}$ for some $q>0$ since the previous cycle operator in $\bz_{\rm red}$ is a power of $\hat A$. This successor is also in $H^0$. Thus $a_{\rm red}-qb_{\rm red}>0$ by (\ref{nn2}), which is a contradiction.
The case  $a=\max\{a,b,c,d\}$ can be treated similarly. \qedohne
\begin{co} \label{Coapp} Let $[B_1,\dots,B_t]$ be a reduced cycle in $H^0$. Then $a_i\not=d_i$ for all indices $i$. Moreover,
$a_i>d_i$ holds if and only if $a_{i+1}<d_{i+1}$. 

If $a_i<d_i$, then $M_i=\hat A^q$ for $$q=\left\{ \begin{array}{ll}
  \left[ \frac{d_i}{c_i}\right],    & \mbox{if } c_i\not=1,\\[1ex]
   d_i-1,  & else,
\end{array}\right.$$
and, if $a_i>d_i$, then $M_i=\hat B^q$ for $$q=\left\{ \begin{array}{ll}
  \left[ \frac{a_i}{b_i}\right],    & \mbox{if } b_i\not=1,\\[1ex]
   a_i-1,  & else.
\end{array}\right.$$
\end{co}
\proof For each $B\in H^0_{\rm red}$, we have $a=\max\{a,b,c,d\}$ and $d=\min\{a,b,c,d\}$ or $a=\min\{a,b,c,d\}$ and $d=\max\{a,b,c,d\}$, see Remark~\ref{Rmax}. Hence $a=d$ would imply $a=b=c=d$, which contradicts $\det B=1$. This proves the first assertion.

Now take an element $B_i$ of the reduced cycle $\bz' =[B_1,\dots,B_t]$ and suppose $a_i<d_i$. By Prop.~\ref{P1red}, the reduced cycle $\bz'$ equals $\bz_{\rm red}$ for a cycle $\bz$ in $H^0$. Consequently, $B_{i+1}=\hat A  ^q B_i\hat A^{-q}\in H^0$ for some $q>0$ and  $B':=\hat A  ^{q+1} B_i\hat A^{-(q+1)}$ is not in $H^0$ by definition of $\bz_{\rm red}$ and Lemma~\ref{LA}. The first condition implies $qc_i<d_i$. We will show that the second condition implies $(q+1)c_i\ge d_i$ and that equality holds only if $c_i=1$, which implies the claimed formula for $q$. Since $a_i+(q+1)c_i>0$ and $c_i>0$, the second condition says that one of the remaining entries of $B'$ is not positive. On the other hand, $\det B'=1$ implies that both remaining entries are non-positive or both are non-negative. Hence both are non-positive, thus  $(q+1)c_i\ge d_i$. If equality holds, then $d_i$ is divisible by $c_i$, hence $\det B'=1$ is divisible by $c_i$, which implies $c_i=1$.\qed
\small\noindent 
Blandine Galiay\\
DER de mathématiques
ENS Paris-Saclay \\
4 avenue des Sciences\\
91190 Gif-sur-Yvette, France. \\
\texttt{blandine.galiay@ens-paris-saclay.fr}\\[2ex]

\small\noindent Ines Kath\\ Institut f\"ur Mathematik und Informatik, Universit\"at Greifswald\\
Walther-Rathenau-Str. 47\\  D-17487 Greifswald, Germany. \\
\texttt{ines.kath@uni-greifswald.de}
\end{document}